\documentclass[11pt]{article}

\usepackage{fullpage,amsfonts,amsmath,amsthm,amssymb,mathrsfs,graphicx,epstopdf}
\usepackage[top=2.54cm, bottom=2.54cm, left=2.54 cm, right=2.54 cm]{geometry}

 \newcommand{\lab}[1]{\label{#1}}                % hides labels

 \usepackage[usenames,dvipsnames]{color}
\newcommand{\fix}[1]{{#1}}

\newcommand{\remove}[1]{}
\newcommand\eqn[1]{(\ref{#1})}

\newcommand{\be}{\begin{equation}}
\newcommand{\bel}[1]{\begin{equation}\lab{#1}\ }
\newcommand{\ee}{\end{equation}}
\newcommand{\bea}{\begin{eqnarray}}
\newcommand{\eea}{\end{eqnarray}}
\newcommand{\bean}{\begin{eqnarray*}}
\newcommand{\eean}{\end{eqnarray*}}

\newtheorem{thm}{Theorem}%[section]
\newtheorem{cor}[thm]{Corollary}

\newtheorem{lemma}[thm]{Lemma}

\newtheorem{claim}[thm]{Claim}
\newtheorem{prop}[thm]{Proposition}

\newtheorem{remark}[thm]{Remark}

\def\proof{\noindent{\bf Proof.~}~}
\def\qed{~~\vrule height8pt width4pt depth0pt}
\def\ss{{\smallskip}}

\newcommand{\ind}[1]{1_{\{#1\}}}
\def\bin{{\bf Bin}}

\def\C{{\mathcal C}}
\def\D{{\mathcal D}}
\def\E{{\mathcal E}}

\def\G{{\mathcal G}}
\def\cH{{\mathcal H}}

\def\P{{\mathcal P}}

\def\ex{{\mathbb E}}
\def\pr{{\mathbb P}}

\def\bfd{{\bf d}}

\def\eps{\epsilon}

\date{}

\title{Subgraph probability of random graphs with specified degrees and applications to chromatic number and connectivity}
\author{Pu Gao\thanks{Research supported by NSERC.} \\
University of Waterloo\\
pu.gao@uwaterloo.ca \and Yuval Ohapkin\\
University of Waterloo\\
yohapkin@uwaterloo.ca
}
\begin{document}
\maketitle

\begin{abstract}
Given a graphical degree sequence ${\bf d}=(d_1,\ldots, d_n)$, let $\G(n,\bfd)$ denote a uniformly random graph on vertex set $[n]$ where vertex $ i$ has degree $d_i$ for every $1\le i\le n$. We give upper and lower bounds on the joint probability of an arbitrary set of edges in $\G(n,\bfd)$, and we link these probability estimates to the corresponding probabilities in the configuration model. Then we show that many existing results of $\G(n,\bfd)$ in the literature can be significantly improved with simpler proofs, by applying this new probabilistic tool. One example we give concerns the chromatic number of $\G(n,\bfd)$. 

In another application, we use these joint probabilities to study the connectivity of $\G(n,\bfd)$. Under some rather mild condition on $\bfd$ --- in particular, if $\Delta^2=o(M)$ where $\Delta$ is the maximum component of $\bfd$ --- we fully characterise the connectivity phase transition of $\G(n,\bfd)$. We also give sufficient conditions for $\G(n,\bfd)$ being connected when $\Delta$ is unrestricted. 
 
\end{abstract}

\section{Introduction}

Given a sequence of nonnegative integers ${\bf d}=(d_1,\ldots, d_n)$, let $\G(n,\bfd)$ denote a uniformly random simple graph on vertex set $[n]$ where vertex $ i$ has degree $d_i$ for every $1\le i\le n$. We also use the same notation $\G(n,\bfd)$ for the support of $\G(n,\bfd)$, i.e.\ the set of graphs with degree sequence $\bfd$. Unless otherwise specified, there is usually no confusion from the context whether $\G(n,\bfd)$ refers to a set or a random graph from this set, and we will specify it if there is confusion. We say $\bfd$ is {\em graphical} if the set $\G(n,\bfd)$ is nonempty. In the special case where $d_i=d$ for every $1\le i\le n$, $\G(n,\bfd)$ is a uniformly random $d$-regular graph.  Random graphs with specified degree sequences are among the most studied random graph models. However, unlike the binomial random graph $\G(n,p)$, edge probabilities in $\G(n,\bfd)$ are correlated. In fact, even estimating the probability of a single edge in $\G(n,\bfd)$ can be challenging. That makes analysing $\G(n,\bfd)$ difficult. A common tool used for analysing $\G(n,\bfd)$ is the configuration model introduced by Bollob\'{a}s~\cite{bollobas1980}. In the configuration model, each vertex is represented as a bin containing $d_i$ points. Take a uniformly random perfect matching over the total $M=\sum_{i=1}^n d_i$ points, and let $G^*$ be the multigraph obtained by taking each pair in the matching as an edge. That is, if points $p$ and $q$ are matched as a pair by the matching then $uv$ is an edge in $G^*$ where $u$ contains $p$ and $v$ contains $q$. $G^*$ is a multigraph because there can be more than one edge between two vertices. A simple counting argument shows that every simple graph with degree sequence $\bfd$ corresponds to the same number of configurations, and thus, $G^*$, conditioned to the set of simple graphs, has the same distribution as $\G(n,\bfd)$. Suppose $\bfd$ is such that the probability that $G^*$ is simple is bounded away from 0 for all large $n$, then any property that holds a.a.s.\ in the configuration model must hold a.a.s.\ in $\G(n,\bfd)$.
 Estimating edge probabilities in the configuration model is quite easy, and hence, many properties of $\G(n,\bfd)$  are obtained by analysing the configuration model. However, the probability that $G^*$ is simple is bounded away from 0 only for $\bfd$ such that $M=\Theta(n)$ and $M_2:=\sum_{i=1}^n d_i^2=O(n)$. This condition significantly restricts the type of results one can get by translating from the configuration model.
 
 The purpose of this work is to develop probabilistic tools for translating configuration model analysis to $\G(n,\bfd)$ analysis for a  rich family of degree sequences. We will illustrate with examples and show how easy it is to improve existing results and prove new properties of $\G(n,\bfd)$ using the new tools we develop.
Let $H$ be a graph on $[n]$. Under some mild conditions on $H$ and $\bfd$, one of our main results shows the following.
 \begin{center}
 The probability that $H$ is a subgraph of $\G(n,\bfd)$ is approximately the probability that all edges in $H$ appear in the configuration model.
 \end{center}
 
 Before formally stating the main results, we define a few necessary terms. 
Given ${\bf d}=(d_1,\ldots, d_n)$, without loss of generality we may assume that $d_1\ge d_2\ge \cdots \ge d_n$. Define
\[
\Delta=d_1,\quad \delta=d_n,\quad M=\sum_{i=1}^n d_i,\quad J({\bf d})=\sum_{i=1}^{\Delta} d_i.
\] 
\fix{Note that for any graph $G$ with degree sequence $\bfd$ and any vertex $v$ in $G$, the number of 2-paths with one end being $v$ is at most $J({\bfd})$.}
For any $S\subseteq [n]$ define
\[
d(S)=\sum_{i\in S} d_i, \quad \Delta_S=\max_{i\in S} d_i.
\]
Given a graph $H$ on vertex set $[n]$, let $e(H)$ be the number of edges in $H$, $\Delta(H)$ be the maximum degree of $H$, and let $\bfd^H=(d^H_1,\ldots, d^H_n)$ denote the degree sequence of $H$.  %For graphs on the same vertex set, e.g.\ on $[n]$, we may treat them as sets of edges. Thus, 
 If two graphs $H_1$ and $H_2$ are both on vertex set $[n]$, we say $H_1$ and $H_2$ are disjoint if their edge sets are disjoint.
Given a graph $H$ on $[n]$, let $H^+$ denote the event that $H \subseteq \G(n,{\bf d})$, i.e.\ $H$ is a subgraph of $\G(n,\bfd)$, and let $H^-$ denote the event that $H\cap \G(n,{\bf d})=\emptyset$, i.e.\ $H$ and $\G(n,{\bf d})$ have disjoint edge sets. Given two degree sequences ${\bf d}$ and ${\bf d}'$ defined on the same set of vertices, we say ${\bf d}\preceq {\bf d}'$ if $d_i\le d'_i$ for every $i$. 

Typically we consider a sequence of degree sequences indexed by $n$, and we are interested in properties of $\G(n,\bfd)$ asymptotically when $n\to\infty$. We say a property holds asymptotically almost surely (a.a.s.) if the probability that the property holds goes to 1 as $n\to \infty$.
We will use standard Landau notation for asymptotic calculations. Given two sequences of real numbers $a_n$ and $b_n$, we say $a_n=O(b_n)$ if there exists a constant $C>0$ such that $|a_n|\le C|b_n|$ for every $n$. We say $a_n=o(b_n)$ if $b_n>0$ eventually (meaning for all sufficiently large $n$) and $\lim_{n\to\infty} a_n/b_n=0$. We write $a_n=\Omega(b_n)$ if $a_n>0$ eventually and $b_n=O(a_n)$. We say $a_n=\omega(b_n)$ if $a_n>0$ eventually and $b_n=o(a_n)$.

\subsection{Main results}

\subsubsection{Subgraph probabilities}
\lab{sec:subgraph_probability}
Our first result concerns the conditional edge probabilities in $\G(n,\bfd)$.

\begin{thm}\label{thm:conditional} 
Let $H_1$ and $H_2$ be two disjoint graphs on $[n]$. Suppose that ${\bf d}^{H_1}\preceq {\bf d}$, and
 $uv\notin H_1\cup H_2$. %Suppose also that \[
%\frac{3J({\bf d})+\Delta(8+d_u^{H_2}+d_v^{H_2})}{M-2e(H_1)}+\frac{2e(H_2)\Delta^2}{(M-2e(H_1))^2}}\le\frac{1}{2}.
%\]
Then,
\begin{eqnarray*}
\pr(uv\in \G(n,{\bf d}) \mid H_1^+, H_2^-)& \le& \frac{(d_u-d^{H_1}_u)(d_v-d^{H_1}_v)}{M-2e(H_1)}\cdot f({\bf d},H_1,H_2)\\
\pr(uv\in \G(n,{\bf d}) \mid H_1^+, H_2^-)&\ge& \frac{(d_u-d^{H_1}_u)(d_v-d^{H_1}_v)}{M-2e(H_1)}\cdot g({\bf d},H_1,H_2) 
\end{eqnarray*}
where 
\begin{eqnarray*}
f({\bf d},H_1,H_2)&=&\left(1-\frac{3J({\bf d})+\Delta(8+d_u^{H_2}+d_v^{H_2})}{M-2e(H_1)}-\frac{2e(H_2)\Delta^2}{(M-2e(H_1))^2}+\frac{(d_v-d_v^{H_1})(d_u-d_u^{H_1})}{M-2e(H_1)}\right)^{-1}\\
g({\bf d},H_1,H_2)&=&\left(1- \frac{2J({\bf d})+6\Delta+2\Delta(H_2)\Delta}{M-2e(H_1)} \right)
\left(1+ \frac{(d_v-d_v^{H_1})(d_u-d_u^{H_1})}{M-2e(H_1)} \right)^{-1},
\end{eqnarray*}
\fix{provided that $f({\bf d},H_1,H_2)>0$.}
\end{thm}

Theorem~\ref{thm:conditional} is proved by a simple switching counting argument. \fix{For the switching argument to work, it is required that there are many options to choose an edge $xy$ where $x$ is not adjacent to some specific vertex $v$. This can be guaranteed if $J({\bfd})$ is relatively small compared to $M$, recalling that $J(\bfd)$ bounds the number of 2-paths starting at $v$. This explains the appearance of terms involving $J(\bfd)/M$ in the above probability bounds.} The proof is presented in Section~\ref{sec:probabilities}.
By setting proper conditions on $\bfd$ we obtain an asymptotic value of $\pr(uv\in \G(n,{\bf d}) \mid H_1^+, H_2^-)$ as follows, which is a direct corollary of Theorem~\ref{thm:conditional}.

\begin{cor}~\lab{cor:conditional}
Let $H_1$ and $H_2$ be two disjoint graphs on $[n]$. Suppose that ${\bf d}^{H_1}\preceq {\bf d}$, and
 $uv\notin H_1\cup H_2$. Suppose further that
 \[
 J(\bfd)+\Delta\cdot\Delta(H_2)=o(M-2e(H_1)),\quad e(H_2)\Delta^2=o\big((M-2e(H_1))^2\big).
 \]
 Then,
 \[
 \pr(uv\in \G(n,{\bf d}) \mid H_1^+, H_2^-)=(1+o(1)) \frac{(d_u-d^{H_1}_u)(d_v-d^{H_1}_v)}{M-2e(H_1)+(d_u-d^{H_1}_u)(d_v-d^{H_1}_v)}.
 \]
\end{cor}

\begin{remark} 
If $\bfd$ is a degree sequence composed of i.i.d.\ copies of a power-law variable with exponent between 2 and 3, then a.a.s.\ $J(\bfd)=o(M)$. In~\cite[Lemma 3]{gao2018}, Van der Hofstad, Southwell, Stegehuis and the first author proved the asymptotic conditional probabilities as in Corollary~\ref{cor:conditional} for such power-law degree sequences when $e(H_1)=O(1)$ and $e(H_2)=0$.  A few other papers studied such conditional probabilities for regular or near-regular degree sequences. Kim, Sudakov and Vu~\cite[Lemma 2.1]{kim2007} obtained the asymptotic value of the conditional probability when $e(H_1)=O(1)$, $e(H_2)=0$ and $d_i\sim d$ for all $i$ where $d=o(n)$ and $d=\omega(1)$. Espuny D\'{i}az, Joos, K\"{u}hn and Osthus~\cite[Lemmas 2.1 and 2.2]{espuny2019} obtained the asymptotic value of such conditional probabilities in a random $d$-regular $r$-uniform hypergraph.  
Our corollary above generalises and strengthens all these results above (compared with~\cite[Lemmas 2.1 and 2.2]{espuny2019} for $r=2$, our corollary has more relaxed conditions on $H_1$ and $H_2$ in addition to relaxed conditions on $\bfd$). Our result is mostly comparable with McKay's enumeration results~\cite{mckay1981}, which yield similar probabilities as in Corollary~\ref{cor:conditional} but require stronger conditions on $\Delta$. These are all for sparse degree sequences. For dense $\bfd$, such conditional probabilities are estimated in~\cite{barvinok2013,isaev2018,mckay2010,gao2020}.
\end{remark}

Repeatedly applying Theorem~\ref{thm:conditional} we obtain the following joint probability bounds for an arbitrary set of edges, under some mild conditions. Given a graph $H$, let $\Lambda(H)=\{v: v\cap x\neq \emptyset\ \mbox{for some $x\in E(H)$}\}$ be the set of vertices that are incident to some edge in $H$. For real number $x$ and nonnegative integer $k$, let $(x)_k=\prod_{i=0}^{k-1}(x-i)$.

\begin{thm}\label{thm:joint}
Assume $H_1$ and $H_2$ are two disjoint graphs on $[n]$.
\begin{enumerate}
\item[(a)] If
$$
\mathfrak{R}({\bf d}, H_1,H_2):=\frac{6J({\bf d})+2\Delta(8+2\Delta(H_2))}{M-2e(H_1)}+ \frac{4e(H_2)\Delta^2}{(M-2e(H_1))^2}\le 1,
$$
then 
\begin{eqnarray*}
\pr\Big(H_1^+\mid H_2^- \Big) \le  \prod_{i=1}^n (d_i)_{d_i^{H_1}}\cdot \prod_{j=1}^{e(H_1)}\frac{1+{\mathfrak{R}}({\bf d}, H_1, H_2)}{ (M-2j+2)}.
%\mathfrak{R}({\bf d}, j, H_2) &=&  1+\frac{6J({\bf d})+2\Delta(8+2\Delta(H_2))}{M-2j+2}+\frac{4e(H_2)\Delta^2}{(M-2j+2)^2}.
\end{eqnarray*}
\item[(b)]
 If
$$
\frak{r}({\bf d},H_1, H_2):=\frac{2J({\bf d})+6\Delta+2\Delta(H_2)\Delta+ \Delta_{\Lambda(H_1)}^2}{M-2e(H_1)} \le 1,
$$
then
\[
\pr\Big(H_1^+\mid H_2^- \Big)\ge \prod_{i=1}^n (d_i)_{d_i^{H_1}}\cdot \prod_{j=1}^{e(H_1)}\frac{1-{\mathfrak{r}}({\bf d}, H_1, H_2)}{ (M-2j+2)}.
\]
\end{enumerate}
\end{thm}
The proof of Theorem~\ref{thm:joint} will be given in Section~\ref{sec:probabilities}.
By setting $H_2=\emptyset$ and imposing conditions on $\bfd$ so that $\frak{R}(\bfd, H_1,\emptyset)=o(1)$, we obtain the following useful upper and lower bounds on the joint probability.

\begin{cor} \label{cor:upperbound}
 Let $H$ be a graph on $[n]$ and assume that ${\bf d}$ is a degree sequence satisfying $J(\bfd)=o(M-2e(H))$. Then
  \begin{equation}
\pr(H^+) \le \prod_{i=1}^n (d_i)_{d_i^{H}} \prod_{i=1}^{e(H)}
\frac{1+o(1)}{M-2i+2}. \lab{graph_prob1}
\end{equation}
If further we have $\Delta^2_{\Lambda(H)}=o(M-2e(H))$ then the above holds with equality.
\end{cor}

\begin{remark}
In the case that~\eqn{graph_prob1} holds with equality, the equality does not give the asymptotic value of $\pr(H^+)$ due to the $1+o(1)$ term in each term of the product. \fix{One can easily get asymptotic value of $\pr(H^+)$ if $e(H)\cdot \mathfrak{R}({\bf d}, H,\emptyset)+e(H)\cdot \mathfrak{r}({\bf d}, H,\emptyset)=o(1)$ by applying Theorem~\ref{thm:joint}. However, lower and upper bounds in Corollary~\ref{cor:upperbound} work for very large $e(H)$, and are sufficient for many applications, as we will show in this paper.}
\end{remark}

\begin{remark}
%Similar bounds on $\pr(H^+)$ are given by McKay in~\cite[Theorem 1.2]{mckay2010} with more restrictive $\bfd$ and $H$. 
An asymptotic estimate of $\pr(H^+)$ was determined in~\cite[Theorem 1.2]{mckay2010} when $\Delta^2 e(H)=o(M)$, which is a stronger condition than $e(H)\cdot \mathfrak{R}({\bf d}, H,\emptyset)+e(H)\cdot \mathfrak{r}({\bf d}, H,\emptyset)=o(1)$. An example of $\bfd$ for which Theorem~\ref{thm:joint} gives an asymptotic estimate for $\pr(H^+)$ provided $e(H)\le n^{\eps}$ for some small constant $\eps>0$, whereas~\cite[Theorem 1.2]{mckay2010} fails even when $e(H)=1$ is the family of i.i.d. power-law degree sequences of exponent $2<\tau<3$ (see more discussions in Section~\ref{sec:examples}).  McKay~\cite[Theorem 2.1]{mckay2010} also determined  $\pr(H^+)$ when $\bfd$ is dense and close to regular, and $H$ is fairly sparse. Our result does not cover these cases.
\end{remark}

\begin{remark} Although probability estimates similar to~\eqn{graph_prob1} appeared in the literature previously~\cite{mckay2010}, such probabilities were not connected before to the corresponding probabilities in the configuration model, and thus were never used to automatically translate the analysis from the configuration model to $\G(n,\bfd)$. This is one of our main contributions of this work.
We compare~\eqn{graph_prob1} with the probability in the configuration model. Let $\sigma^*$ denote a uniformly random perfect matching over the $M$ points produced by the configuration model, and let $G^*$ be the multigraph corresponding to $\sigma^*$. Given a graph $H$ on $[n]$ with $\bfd^H\preceq \bfd$, let $\P(H)$ be the set of matchings $\sigma$ of size $e(H)$ over the set of $M$ points whose corresponding graph is $H$. Then $|\P(H)|=\prod_{i=1}^n (d_i)_{d_i^{H}}$ and $\pr(\sigma\subseteq \sigma^*)=\prod_{i=1}^{e(H)}
\frac{1}{M-2i+1}$. Hence, 
\begin{equation}
\pr(H\subseteq G^*)\le \sum_{\sigma\in \P(H)}\pr\big(\sigma\subseteq \sigma^*\big)= \prod_{i=1}^n (d_i)_{d_i^{H}} \prod_{i=1}^{e(H)}
\frac{1}{M-2i+1}.\lab{CM_prob1}
\end{equation}
Thus, under the condition $J(\bfd)=o(M-2e(H))$, our upper bound in~\eqn{graph_prob1} differs from the corresponding probability bound~\eqn{CM_prob1} in the configuration model by a relative $1+o(1)$ factor in each term of the product. Such an approximation is enough to translate a lot of configuration model analysis to $\G(n,\bfd)$. See examples in Sections~\ref{secsec:chrom}.
\end{remark}

Another advantage of using the configuration model is that, it is very easy to bound the probability of having a certain number of edges joining two sets of vertices. 
Given two subsets of vertices $S_1,S_2\subseteq [n]$, let $e(S_1, S_2)$ be the number of edges with one end in $S_1$ and the other end in $S_2$. When $S_1=S_2$, $e(S_1,S_2)$ is simply the number of edges induced by $S_1$. In the configuration model, the probability that $e(S_1,S_2)\ge \ell$ is at most 
\begin{equation}
\binom{d(S_1)}{\ell} \frac{(d(S_2))_{\ell}}{\prod_{i=1}^\ell (M-2i+1)}. \lab{CM_prob2}
\end{equation}
In the following corollary we show a similar bound for this probability in $\G(n,\bfd)$.
\begin{cor} \label{cor:join}
Suppose $S_1, S_2 \subseteq [n]$. Let $1\le \ell<M/2$ be an integer.  Assume that ${\bf d}$ satisfies $J(\bfd) = o(M-2\ell)$. Then,%the probability that $e(S_1 , S_2) \ge \ell$ in $\G(n, {\bf d})$ is at most 
\begin{equation}
\pr(e(S_1 , S_2) \ge \ell)\le\binom{d(S_1)}{\ell} (d(S_2))_{\ell} \left(\prod_{i=1}^{\ell} \frac{1+o(1)}{M-2i+2}\right) \le\binom{d(S_1)}{\ell} \frac{(d(S_2))_{\ell}}{(M/2)_{\ell} (2+o(1))^{\ell}}. \lab{graph_prob2}%\le \left((e/2+o(1))\frac{d(S_1)d(S_2)}{\ell(M/2-\ell)}\right)^{\ell}.
\end{equation}
%$\left(\frac{O(D(S_1)D(S_2))}{\ell(M-2\ell)}\right)^{\ell}$. Moreover, if $\min\{D(S_1), D(S_2)\} \le M-\ell+1$ (in particular, when $S_1$ and $S_2$ are disjoint) this bound can be improved to $ \binom{D(S_1)}{\ell} \frac{((1+o(1))D(S_2))^{\ell}}{(M)_{\ell}}$.
\end{cor}

\proof Let $\mathcal{H}_{\ell}$ denote the set of graphs on $[n]$ having exactly $\ell$ edges, all of which have exactly one end in $S_1$ and the other end in $S_2$. Let $\{i_1 ,\cdots , i_k\}$ denote the set of vertices in $S_1\cup S_2$. By the union bound and Corollary~\ref{cor:upperbound}, %
%the probability that $H \subseteq \mathcal{G}(n, {\bf d})$ for some $H \in \H_{\ell}$ is at most 
\[
\pr(e(S_1 , S_2) \ge \ell) \le \sum_{H\in\cH_{\ell}} \pr(H\subseteq \G(n,{\bf d}))\le
\left(\prod_{i=1}^{\ell} \frac{1+o(1)}{M-2i+2}\right) \sum_{H \in \cH_{\ell}} \prod_{j=1}^k (d_{i_j})_{d_{i_j}^H}.% \le \left(\frac{1+o(1)}{M-2\ell}\right)^{\ell} \sum_{H \in \H_{\ell}} \prod_{j=1}^k (d_{i_j})_{d_{i_j}^H}
\]
Next we give a combinatorial interpretation of $\sum_{H \in \mathcal{H}_{\ell}} \prod_{j=1}^{k} (d_{i_j})_{d_{i_j}^{H}}$ above. Represent each vertex $i_j$ in $S_1\cup S_2$ by a bin containing $d_{i_j}$ points. Given $H\in\cH_{\ell}$, how many matchings of size $\ell$ over $2\ell$ out of the total $d(S_1\cup S_2)$ points are there so that if $uv$ is an edge in $H$ then there is a point $p$ in bin $u$ and a point $q$ in bin $v$ such that $p$ and $q$ are matched by the matching? It is easy to see that there are exactly $\prod_{j=1}^{k} (d_{i_j})_{d_{i_j}^{H}}$ such matchings. Hence, $\sum_{H \in \mathcal{H}_{\ell}} \prod_{j=1}^{k} (d_{i_j})_{d_{i_j}^{H}}$ is bounded above by the total number of size-$\ell$ matchings where every pair $(p,q)$ in the matching is of the form that $p$ is in some bin in $S_1$ and $q$ is in some bin in $S_2$. There are $\binom{d(S_1)}{\ell}$ ways to choose the $\ell$ ends that are in bins in $S_1$ and there are $(d(S_2))_{\ell}$ ways to choose the other ends from bins in $S_2$ and match them to the $\ell$ ends chosen before. Hence,
\[
\sum_{H \in \mathcal{H}_{\ell}} \prod_{j=1}^{k} (d_{i_j})_{d_{i_j}^{H}} \le \binom{d(S_1)}{\ell} (d(S_2))_{\ell},
\]
and thus,
\[
\pr(e(S_1 , S_2) \ge \ell) \le
\left(\prod_{i=1}^{\ell} \frac{1+o(1)}{M-2i+2}\right) \binom{d(S_1)}{\ell} (d(S_2))_{\ell}=\binom{d(S_1)}{\ell}\frac{(d(S_2))_{\ell}}{(M/2)_{\ell} (2+o(1))^{\ell}},
\]
as $\prod_{i=1}^{\ell}(M-2i+2)=(M/2)_{\ell} 2^{\ell}$.  \qed
\smallskip

\subsubsection{Chromatic number of $\G(n,\bfd)$}
\lab{secsec:chrom}
Let $\chi(G)$ denote the chromatic number of graph $G$, i.e.\ the minimum number of colours required to colour vertices of $G$ so that all pairs of adjacent vertices receive distinct colours.
It is known that a.a.s.\ the chromatic number of a random $d$-regular graph is asymptotically $d/2\ln d$, for $d=\omega(1)$ and $d=o(n)$; see~\cite{frieze1992, cooper2002random,krivelevich2001}. In the paper~\cite{frieze2007} Frieze, Krivelevich and Smyth asked under what conditions on $\bfd$ would we have a.a.s.\ $\chi(\G(n,\bfd))=\Theta(d/\ln d)$, where $d=M/n$ is the average degree of graphs in $\G(n,\bfd)$. Let $D_k=\sum_{i=1}^k d_i$ and $M_2=\sum_{i=1}^n d_i^2 $. It was shown in~\cite{frieze2007} that if
\begin{itemize}
\item[(A1)] there exist constants $1/2<\alpha<1$, $\eps, K_0>0$ such that $D_k\le K_0 dn (k/n)^{\alpha}$ for all $1\le k\le \eps n$;
\item[(A2)] $\Delta^5=o(M_2)$,
\end{itemize}
then a.a.s.\ $\chi(\G(n,\bfd))=O(d/\ln d)$.
On the other hand, if
\begin{itemize}
\item[(A3)] $\Delta^4=o(M)$,
\end{itemize}
then a.a.s.\ $\chi(\G(n,\bfd))=\Omega(d/\ln d)$.
\smallskip

We significantly relax the conditions (A2) and (A3) and obtain the following result for $\chi(\G(n,\bfd))$.
\begin{thm}\lab{thm:chrom}
If $\bfd$ satisfies condition (A1) and $\Delta=o(n)$ then a.a.s.\ $\chi(\G(n,\bfd))=O(d/\ln d)$. If $J(\bfd)=o(M)$ then a.a.s.\ $\chi(\G(n,\bfd))=\Omega(d/\ln d)$.
\end{thm}
The proof of Theorem~\ref{thm:chrom} will be given in Section~\ref{sec:chrom}, which is obtained by translating the existing analysis~\cite{frieze2007} from the configuration model to $\G(n,\bfd)$. We believe that many other results of $\G(n,\bfd)$ can be obtained or improved in a similar manner. For instance, %the giant component threshold for $\G(n,\bfd)$ for the family of the so-called ``well-behaved'' $\bfd$ was determined by Molloy and Reed~\cite{molloy1995}. This result was recently extended to all degree sequences by  Joos, Perarnau, Rautenbach and Reed~\cite{joos2018} and the proof uses the switching method on the random graphs. 
the order of the largest component of $\G(n,\bfd)$ was determined by Molloy and Reed~\cite{molloy1998}  for the so-called ``well-behaved'' degree sequences. Their proof relies on an analysis in the configuration model.  We believe that most of the analysis can be immediately translated to $\G(n,\bfd)$, \fix{for a much larger family of $\bfd$,} by using the conditional edge probabilities in Theorem~\ref{thm:conditional}. \fix{Theorem~\ref{thm:conditional} requires certain conditions on $\bfd$. Joos et. al.~\cite{joos2018} extended the work~\cite{molloy1998} to every feasible degree sequence $\bfd$, by using extensive switching techniques; some of which are similar to the ones we use in the paper}.  We also believe that the new probabilistic tools developed in Section~\ref{sec:subgraph_probability} will be useful in studying other properties of $\G(n,\bfd)$. We give another example in Section~\ref{secsec:connectivity}.

\subsubsection{Connectivity transition of $\G(n,\bfd)$}
\lab{secsec:connectivity}

The connectivity is one of the best studied graph properties for random graphs.  Erd{\H{o}}s and R\'{e}nyi~\cite{erdHos1964} determined the threshold of the connectedness for $\G(n,p)$. Indeed, for every fixed integer $k\ge 1$,  Erd{\H{o}}s and R\'{e}nyi ~\cite{erdHos1964} determined when $\G(n,p)$ becomes $k$-connected.  Random graph $\G(n,p)$ becomes connected when isolated vertices disappear, which happens when $p\approx \log n/n$. For $p$ in this range, the average degree of the random graphs is around $\log n$ and there are very few vertices of degree one or two, and these vertices are pair-wise far away in graph distance. Consequently the vertices of degree one or two do not affect the connectedness of $\G(n,p)$. The most natural sparser random graph model for the study of connectivity would be $\G(n,\bfd)$, and it is natural to ask when are such random graphs connected. As we will see, the vertices of degree one or two play crucial roles for the connectedness of $\G(n,\bfd)$.

The first work on the connectivity of $\G(n,\bfd)$ was by Wormald~\cite{wormald1981} in 1981. In this pioneering work, the author studied the connectivity of $\G(n,\bfd)$ where $\delta\ge 3$ and $\Delta$ is bounded by some absolute constant $R$ (i.e.\ $R$ does not depend on $n$), and proved that a.a.s.\ the connectivity of $\G(n,\bfd)$ is equal to $\delta$ for such degree sequences.  Frieze~\cite{frieze1988} studied the connectivity of random $d$-regular graphs where $3\le d =o(n^{0.2})$, and proved that a.a.s.\ a random $d$-regular graph is $d$-connected, for $d$ in the aforementioned range. 
Later, Cooper, Frieze and Reed~\cite{cooper2002} extended this result to $3\le d\le cn$ where $c>0$ is a sufficiently small constant. {\L}uczak~\cite{luczak1989} extended their work to non-regular degree sequences, and considered also degree sequences permitting $\delta=2$. {\L}uczak showed that, for any $\bfd$ where $\delta\ge 3$ and $\Delta\le n^{0.01}$ then a.a.s.\ $\G(n,\bfd)$ is $\delta$-connected. When $\delta=2$ and $\Delta\le n^{0.01}$ he characterised the structure of $\G(n,\bfd)$ and determined when is $\G(n,\bfd)$ a.a.s.\ 2-connected.
In a more recent work, Federico and Van der Hofstad~\cite{federico2016} considered degree sequences permitting $\delta=1$ and fully charactersized the connectivity transition of $\G(n,\bfd)$ for $\bfd\in {\D}$, where ${\D}=\{\bfd: M=\Theta(n), \sum_{i=1}^n d_i^2=O(n)\}$. Let $n_1$ be the number of \fix{entries} in $\bfd$ with value 1, and $n_2$ the number of \fix{entries} in $\bfd$ with value 2.  Federico and Van der Hofstad showed that for $\bfd\in {\D}$ that satisfies some additional ``smoothness'' condition, $\G(n,\bfd)$ is a.a.s.\ connected if $n_1=o(\sqrt{n})$ and $n_2=o(n)$, and $\G(n,\bfd)$ is disconnected if $n_1=\omega(\sqrt{n})$. All the work that we have discussed so far are for $\bfd$ where either the maximum degree is not large (at most $n^{0.01}$), or $\bfd$ corresponds to a regular degree sequence, and the degree is nearly sublinear (at most $cn$ for some sufficiently small $c$). For $d$ linear in $n$,  Krivelevich, Sudakov, Vu and Wormald~\cite{krivelevich2001} proved several properties of random $d$-regular graphs, including the connectivity. Recently, Isaev, McKay and the first author~\cite{gao2020} proved several properties including the connectivity of $\G(n,\bfd)$ for near-regular $\bfd$ where $d=\omega(\log n)$ and  $d_i\sim d$ for every $i$.  

In this work, we characterise the connectivity transition of $\G(n,\bfd)$ for a much larger family of degree sequences. For the family of $\bfd$ where $J(\bfd)=o(M)$ (in particular, when $\Delta^2=o(M)$) we fully characterise the phase transition of the connectedness of $\G(n,\bfd)$. When $\Delta$ is unrestricted we give sufficient conditions under which $\G(n,\bfd)$ is a.a.s.\ connected. 

 We only consider degree sequences where $\delta\ge 1$ since otherwise $\G(n,\bfd)$ is disconnected trivially.
Given the degree sequence $\bfd$ where $\Delta=d_1\ge d_2\ge \cdots \ge d_n=\delta\ge 1$, define
\begin{eqnarray*}
&&n_1= \sum_{i=1}^n \ind{d_i=1},\quad n_2=\sum_{i=1}^n \ind{d_i=2}.
\end{eqnarray*}
%Given $S\subseteq [n]$, let ${\overline S}=[n]\setminus S$. Let $e(S)$ denote the number of edges induced by $S$ and $e(S,{\overline S})$ denote the number of edges with exactly one end in $S$ and the other end in ${\overline S}$. For a set of edges $H$ we use $e(H)$ to denote $|H|$, i.e.\ the number of edges in $H$. 

\begin{thm}\lab{thm:small}
Assume ${\bfd}$ is a graphical degree sequence with $J(\bfd)=o(M)$. Let $c>0$ be a fixed constant.
\begin{enumerate}
\item[(a)] If $n_1=o(\sqrt{M})$ and $n_2=o(M)$ then a.a.s.\ $\G(n,\bfd)$ is connected.
\item[(b)] If $n_1=\omega(\sqrt{M})$ then a.a.s.\ $\G(n,\bfd)$ is disconnected.
\item[(c)] If $n_1\ge c\sqrt{M}$ or $n_2\ge c M$ then there exists $\eta=\eta(c)>0$ such that for all sufficiently large $n$,
\[
\pr(\G(n,\bfd)\ \text{disconnected}) \ge \eta.
\]
\end{enumerate}
\end{thm}

Since $J(\bfd)\le \Delta^2$, we immediately have the following corollary.
\begin{cor}
Corollaries~\ref{cor:upperbound} and~\ref{cor:join}, and Theorem~\ref{thm:small} hold if $J(\bfd)=o(M)$ is replaced by $\Delta^2=o(M)$.
\end{cor}

Next, we deal with degree sequences where $\Delta$ is rather large.
Define 
\[
\cH=\{i: d_i\ge \sqrt{M}/\log M\}.%,\quad {\widetilde M}=M-2d(\cH).
\]
Our next result gives sufficient conditions for the connectedness of $\G(n,\bfd)$.
\begin{thm}\lab{thm:H} Assume $M-2d(\cH)=\Omega(M)$.
If $n_1=o(\sqrt{M})$ and $n_2=o(M)$ then a.a.s.\ $\G(n,\bfd)$ is connected.
\end{thm}

The proofs for Theorems~\ref{thm:small} and~\ref{thm:H} will be presented in Section~\ref{sec:connectivity}. We note here that
the conditions in Theorem~\ref{thm:H} are not necessary. 
 We can easily make up $\bfd$ where $d_1=n-1$, $M=\Theta(n)$ and a linear number of vertices have degree 1 (or of degree 2).  
Conditions in Theorem~\ref{thm:H} are not satisfied but $\G(n,\bfd)$ is always connected for such degree sequences.

\subsection{Degree sequences satisfying $J(\bfd)=o(M)$}
\label{sec:examples}

%The family of degree sequences covered by Theorems~\ref{thm:small} and~\ref{thm:H} is much richer than those in the existing work~\cite{wormald1981,frieze1988, luczak1989, cooper2002, federico2016}. 
A rich family of degree sequences satisfies the condition $J(\bfd)=o(M)$ for which we may apply the probabilities in Theorems~\ref{thm:conditional} and~\ref{thm:joint}. For instance, it covers all regular sublinear degree sequences, i.e.\ $d_i=d$ for all $i$ and $d=o(n)$.
We give two additional examples below that might be of interest in applications.
\begin{itemize}
\item $\Delta=o(n)$ and a linear (in $n$) number  of vertices have degree $\Omega(\Delta)$.
\item $\bfd$ is composed of i.i.d.\ power-law variables of exponent $\tau>2$, conditioned to even sum.
\end{itemize}
%The next type of degree sequences satisfy $M'=\Omega(M)$.
%\begin{itemize}
%\item  There exists constant $\eps>0$ with $(2+\eps)a<\bar d$, where $\bar d$ denote $\norm{\bfd}_1/n$ and $a$ be the number of vertices of degree at least $\sqrt{\bar d n}/\log n$. 
%\end{itemize}
The reader may wonder what types of degree sequences do not satisfy $J(\bfd)=o(M)$. Certainly, regular degree sequences with linear degrees do not satisfy this condition. More generally,
if there is a linear (in $n$) number  of vertices with degree $\Theta(n)$, then that degree sequence does not satisfy $J(\bfd)=o(M)$. \fix{There are certainly also examples where the maximum degree is sublinear; for instance, consider $\bfd$ where $\Delta$ vertices have degree $\Delta$, where $\Delta^2=\Omega(n)$, and the remaining vertices all have degree $O(1)$.}
\smallskip

We will 
prove Theorems~\ref{thm:conditional} and~\ref{thm:joint} in Section~\ref{sec:probabilities}. The proof of Theorem~\ref{thm:chrom} will be given in Section~\ref{sec:chrom} and the proofs for Theorems~\ref{thm:small} and~\ref{thm:H} will be presented in Section~\ref{sec:connectivity}.

\section{Proof of Theorem~\ref{thm:chrom}: chromatic number}
\lab{sec:chrom}

We first briefly sketch the proof in~\cite{frieze2007}. For the upper bound, the authors first obtained an $O(d/\ln d)$ upper bound for the multigraph $G^*$ from the configuration model. A multigraph is properly coloured if every pair of adjacent and  distinct vertices are coloured differently. Only Condition (A1) is needed for this part of the proof. Then they applied a sequence of switching operations which repeatedly switch away the loops and multiple edges in $G^*$. Then they proved that (a), every simple graph is obtained with asymptotically the same probability after applying the switchings; (b) if $H$ is the graph induced the by set of edges added during the switchings, then a.a.s.\ $\chi(H)=O(1)$. It follows immediately that a.a.s.\ the chromatic number of $\G(n,\bfd)$ is at most $\chi(G^*)\cdot \chi(H)=O(d/\ln d)$. Condition (A2) is needed to guarantee (a) and (b).

For the lower bound, they proved that for any partition $\sigma$ of vertices into $t=b\cdot d/\ln d$ parts, where $b>0$ is a sufficiently small constant, the probability that $\sigma$ specifies a proper $t$-colouring of $G^*$ is at most $t^{-2n}$. Condition (A3) was applied to obtain a lower bound on the probability that $G^*$ is a simple graph. When $\Delta^2=o(n)$, the probability $t^{-2n}$ is small enough to beat the union bound over all such partitions $\sigma$, and the inverse of the probability that $G^*$ is simple.  Note that (A3) implies $\Delta^2=o(n)$.

To prove Theorem~\ref{thm:chrom}, we carry all analysis from~\cite{frieze2007} for $G^*$ to $\G(n,\bfd)$. Both the upper and lower bound proofs for $\chi(G^*)$ in~\cite{frieze2007} use upper bounds on the probability of $G^*$ containing some set of edges $H$ where $M-2e(H)=\Omega(M)$. The upper bound of $\chi(G^*)$ follows by~\cite[Lemmas 3.1--3.3]{frieze2007}. \fix{\cite[Lemmas 3.1 and 3.3]{frieze2007}} hold for $\G(n,\bfd)$ with exactly the same proofs, by applying inequalities~\eqn{graph_prob1} and~\eqn{graph_prob2} instead of~\eqn{CM_prob1} and~\eqn{CM_prob2}. The additional $1+o(1)$ factors in~\eqn{graph_prob1} and~\eqn{graph_prob2}, compared with~\eqn{CM_prob1} and~\eqn{CM_prob2}, do not affect the proof (in fact, any constant factor would do). As no switching analysis is required any more, we do not need condition (A2). Instead, we need $J(\bfd)=o(M)$ in order to apply~\eqn{graph_prob1} and~\eqn{graph_prob2}. This is guaranteed by our assumption (A1) and $\Delta=o(n)$ as follows:  by (A1), $J(\bfd)/M=O((\Delta/n)^{\alpha})$ which is $o(1)$ if $\Delta=o(n)$. \fix{The proof for~\cite[Lemma 3.2]{frieze2007} uses Azuma-Hoeffding type inequalities which do not hold in general for $\G(n,\bfd)$. In Section~\ref{sec:modification} we describe how to modify the concentration argument in ~\cite[Lemma 3.2]{frieze2007}.}

The same translation of analysis holds for the lower bound proof. As we are working on $\G(n,\bfd)$ instead of $G^*$, it is sufficient if the probability $t^{-2n}$ beats the union bound over the total number of such partitions. This is always the case as there can be at most $t^n$ partitions into $t$ parts. Hence, for the lower bound we do not need conditions (A3) or $\Delta^2=o(n)$  any more. Instead, we impose $J(\bfd)=o(M)$ which validates the application of the probability bounds~\eqn{graph_prob1} and~\eqn{graph_prob2}.  

\fix{
\subsection{Modification of~\cite[Lemma 3.2]{frieze2007}}
\label{sec:modification}
We refer the readers to~\cite{frieze2007} for the definitions of $\eps_1$, $t_0$, $t_1$, $\Delta_t$, subsets $I_0,\ldots, I_{t_1}$ of $[n]$. Notice that the degrees are in non-decreasing order in~\cite{frieze2007}. Let $G_t$ be the subgraph induced by $I_t$, $\tilde d_k$ the degree of vertex $k\in I_t$ in $G_t$, and let $B_t=\{k\in I_t: \tilde d_k\ge \Delta_t\}$ for every $0\le t\le t_1$. Let $Z_t=|B_t|$ and let $Z=Z_0+\ldots +Z_{t_1}$. \cite[Lemma 3.2]{frieze2007} states that a.a.s.\ $Z\le \eps_1 n/2$. Most of the proofs carry directly to $\G(n,\bfd)$, e.g.\ $B_t=\emptyset$ for all $0\le t\le t_0$. For $t>t_0$, $\tilde d_k$ is stochastically dominated by the binomial random variable $\bin\left(d_k,(1+o(1))\frac{D_{n/3^t}}{dn}\right)$ (instead of $\bin\left(d_k, \frac{D_{n/3^t}}{dn}\right)$ as in~\cite{frieze2007}). Then, with exactly the same proof and the same choice of $K_2'$, $\pr(\tilde d_k\ge \Delta_t)\le e^{-\Delta_t}$. It follows immediately then that
\[
\ex Z \le n \sum_{t=1}^{t_1} 3^{-t} e^{-\Delta_t} < n e^{-\Delta_{t_1}}.
\]
By the various choices of parameters such as $\eps_1$, $t_1$, $\Delta_t$, etc.\ it was proved that the above is at most $\eps_1 n/4$.
Then, an Azuma-Hoeffding type of inequality applies to show that $Z\le 2\ex Z$ a.a.s.\ for $G^*$. In this section, we show $Z\le 2\ex Z$ a.a.s.\ using a different concentration argument.

First notice that $D_{n/3^t}/dn \le K_0 3^{-\alpha t}$ which implies that
\[
\frac{D_{n/3^t}}{dn- 2 D_{n/3^t}} < 2K_0 3^{-\alpha t},
\]
as $t>t_0=\lfloor \log_3 1/\eps \rfloor$ and $\eps>0$ can be assumed sufficiently small without loss of generality. Now, expose $\tilde d_k$ one by one for each $k\in I_t$. Conditional on  $(\tilde d_h:\ h\in S)$ for any $S\subseteq I_t$, $\tilde d_k$ is stochastically dominated by the binomial random variable $\bin(d_k,(1+o(1))\frac{D_{n/3^t}}{dn-2D_{n/3^t}})$ by Theorem~\ref{thm:conditional}. By the same calculations as in~\cite{frieze2007} with $K_2'=4K_0^2 e^2 3 ^{1-\alpha}$ (instead of $K_0^2 e^2 3 ^{1-\alpha}$ as in~\cite{frieze2007}), we have 
\[
\pr(\tilde d_k \ge \Delta_t \mid (\tilde d_h:\ h\in S)) \le e^{-\Delta_t},\quad \mbox{for any $S\subseteq I_t$ where $k\notin S$.}
\]
Hence, $Z_t$ is stochastically dominated by $\bin(|I_t|, e^{-\Delta_t})$. It follows immediately by the Chernoff bound that a.a.s.\ $Z_t\le 1.5 \ex Z_t \cdot\ind{\ex Z_t\ge \sqrt{n}} + n^{3/4} \cdot \ind{\ex Z_t <\sqrt{n}}$ for every $t\le t_1$.
It follows then that a.a.s.\
$Z=\sum_{t=0}^{t_1} Z_t \le 1.5 \ex Z + n^{3/4}t_1 < 2\ex Z$, as $t_1=O(\log n)$ and $\ex Z=\Theta(n)$. \qed

}

\section{Proofs of Theorems~\ref{thm:conditional} and~\ref{thm:joint}}
\label{sec:probabilities}

We first prove Theorem~\ref{thm:joint} assuming Theorem~\ref{thm:conditional}. \medskip

{\em Proof of Theorem~\ref{thm:joint}. } Let $e_1=u_1v_1,\ldots, e_h=u_hv_h$ be an enumeration of the set of edges in $H_1$ where $h=e(H_1)$. Let $G_0=\emptyset$ and $G_j =G_{j-1}\cup\{e_j\}$ for every $1\le j\le h$. Then, 
\[
\pr\Big(H_1^+\mid H_2^- \Big)=\prod_{j=1}^h \pr(e_j\mid G_{j-1}^+,H_2^-).
\]
By Theorem~\ref{thm:conditional},
\begin{eqnarray*}
\prod_{j=1}^h \pr(e_j\mid G_{j-1}^+,H_2^-) &\le& \prod_{j=1}^h (d_{u_j}-d_{u_j}^{G_{j-1}})(d_{v_j}-d_{v_j}^{G_{j-1}}) \frac{f({\bf d}, G_{j-1}, H_2)}{M-2j+2}\\
&=&\prod_{i=1}^n (d_i)_{d_i^{H_1}}\cdot \prod_{j=1}^{h}\frac{f({\bf d}, G_{j-1}, H_2)}{ (M-2j+2)},
\end{eqnarray*}
%The contribution of vertex $i \in [n]$ in the product $\prod_{j=1}^h (d_{u_j}-d_{u_j}^{H_{j-1}^-})(d_{v_j}-d_{v_j}^{H_{j-1}^-})$ is 
%\[ 
%(d_i)_{d_i^{H^+}} = d_i(d_i-1)\cdots (d_i-d_i^{H^+}+1),
%\]
%and
where
\begin{eqnarray*}
f({\bf d}, G_{j-1},H_2) &\le& \left(1-\frac{3J({\bf d})+\Delta(8+2\Delta(H_2))}{M-2j+2}-\frac{2e(H_2)\Delta^2}{(M-2j+2)^2}\right)^{-1}\\ 
&\le&1+\frac{6J({\bf d})+2\Delta(8+2\Delta(H_2))}{M-2j+2}+\frac{4e(H_2)\Delta^2}{(M-2j+2)^2}\le 1+\mathfrak{R}({\bf d}, H_1, H_2), 
\end{eqnarray*}
for every $1\le j\le h$. The second inequality above holds by the fact that $(1-x)^{-1}\le 1+2x$ for all $x\in[0,1/2]$ and the assumption that $\mathfrak{R}({\bf d}, H_1, H_2)\le 1$. This yields our upper bound on $\pr(H_1^+\mid H_2^-)$.
Again by the lower bound in Theorem~\ref{thm:conditional},
\begin{eqnarray*}
\prod_{j=1}^h \pr(e_j\mid G_{j-1}^+,H_2^-) &\ge& \prod_{j=1}^h (d_{u_j}-d_{u_j}^{G_{j-1}})(d_{v_j}-d_{v_j}^{G_{j-1}}) \frac{g({\bf d}, G_{j-1}, H_2)}{M-2j+2}\\
&=&\prod_{i=1}^n (d_i)_{d_i^{H_1}}\cdot \prod_{j=1}^{h}\frac{g({\bf d}, G_{j-1}, H_2)}{ (M-2j+2)},
\end{eqnarray*}
where 
\begin{eqnarray*}
g({\bf d}, G_{j-1},H_2) &=& \left(1- \frac{2J({\bf d})+6\Delta+2\Delta(H_2)\Delta}{M-2j+2} \right)
\left(1+ \frac{(d_{v_j}-d_{v_j}^{G_{j-1}})(d_{u_j}-d_{u_j}^{G_{j-1}})}{M-2j+2} \right)^{-1}\\
&\ge & 1- \frac{2J({\bf d})+6\Delta+2\Delta(H_2)\Delta+ d_{v_j}d_{u_j}}{M-2j+2} \ge 1-\frak{r}({\bf d},H_1, H_2),
\end{eqnarray*}
for every $1\le j\le h$. The second inequality above holds by the fact that $(1+x)^{-1}\ge 1-x$ for all $x\ge 0$. This yields our lower bound on $\pr(H_1^+\mid H_2^-)$.
\qed

%Using the inequality $1+x \le (1-x)^{-1}$ for $0 \le x < 1$ and the assumption that $$6J({\bf d})+2\Delta(8+2\Delta(H^-)) < M-2e(H^+),$$ this implies
%\[
%\prod_{j=1}^h \frac{f({\bf d},H_{j-1}^+,H^-)}{M-2j+2} \le \frac{1}{\prod_{j=1}^{e(H^+)} (M-2j+2-6J({\bf d})-2\Delta(8+2\Delta(H^-)))} \left( 1+\frac{4e(H^-)\Delta^2}{(M-2e(H^+))^2} \right)^{e(H^+)},
%\]
%so
%\[
%\pr\Big(A(H^+),{\overline A}(H^-) \Big) \le \pr({\overline A}(H^-)) \frac{\prod_{i=1}^n (d_i)_{d_i^{H^+}}}{\prod_{i=1}^{e(H^+)} (M-2i-6J({\bf d})-2\Delta(8+2\Delta(H^-)))} \left( 1+\frac{4e(H^-)\Delta^2}{(M-2e(H^+))^2} \right)^{e(H^+)}
%\]

\medskip

{\em Proof of Theorem~\ref{thm:conditional}. } Let $\G$ denote the set of graphs $G$ on $[n]$ with degree sequence $d_1, \dots, d_n$, such that $H_1 \subseteq G$, $G \cap H_2 = \emptyset$. Let $S$ denote the set of graphs in $\G$ that contain $uv$ as an edge, and let $\overline{S}=\G\setminus S$. Then,

$$
\pr\big(uv \in \G(n, {\bf d}) \mid H_1^+, H_2^-\big) = \frac{|S|}{|S|+|\overline{S}|} = \frac{1}{1+|\overline{S}|/|S|}
$$

We will obtain upper and lower bounds on the ratio $|\overline{S}|/|S|$ by analysing {\em switchings} that relate graphs in $S$ to graphs in $\overline{S}$. We first define the switching.

\begin{figure}
  \[
  \includegraphics[scale=1]{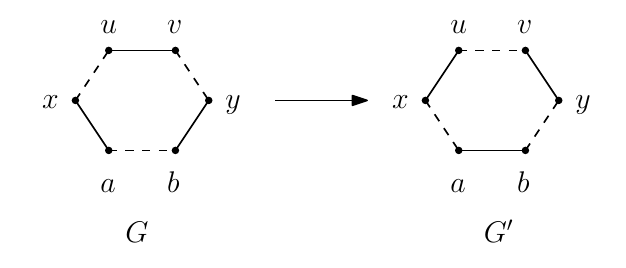}  
  \]
  \caption{Forward switching}
  \label{f:switch}
\end{figure}

Given $G \in S$, a {\em forward switching} specifies an ordered $4$-tuple $(x,a,y,b) \in [n]^4$ satisfying the following conditions:

\begin{itemize}
    \item[(1)] $u, v, x, y, a, b$ are all distinct, except $x = y$ is permitted.
    \item[(2)] $xa$ and $yb$ are edges in $G \setminus H_1$.
    \item[(3)] None of $xu$, $yv$, and $ab$ are edges in $G \cup H_2$. 
\end{itemize}
Then the forward switching converts $G$ to a graph $G' \in \overline{S}$ by deleting the edges $uv$, $xa$ and $yb$ from $G$ and adding the edges $xu$, $yv$ and $ab$. See Figure~\ref{f:switch} for an illustration, where solid lines denote edges in the graph and dashed lines denote non-edges.

Let $f(G)$ denote the number of forward switchings that can be applied to $G$. We will show the following upper and lower bounds on $f(G)$:
\begin{claim}\label{c:forward}
\begin{eqnarray*}
(a)&&f(G)\le (M-2e(H_1))^2\\
(b)&&f(G)\ge  (M-2e(H_1))^2
\left(1-\frac{3J({\bf d})+\Delta(8+d_u^{H_2}+d_v^{H_2})}{M-2e(H_1)}-\frac{2e(H_2)\Delta^2}{(M-2e(H_1))^2}\right).
\end{eqnarray*}
\end{claim}

Next, given $G'\in \overline{S}$, we count the number of forward switchings that can produce $G'$. In order to do so, we define a \textit{backward switching} on $G'$ as an ordered $4$-tuple $(x,a,y,b) \in [n]^4$ satisfying the following:
\begin{itemize}
    \item[(1')] $u, v, x, y, a, b$ are all distinct, except $x = y$ is permitted.
    \item[(2')] $xa$ and $yb$ are not edges in $G' \cup H_2$.
    \item[(3')] $xu$, $yv$, and $ab$ are edges in $G' \setminus H_1$.
\end{itemize}
 Then the backward switching  deletes the edges $xu$, $yv$, and $ab$, and adds the edges $uv$, $xa$ and $yb$. 

Obviously, a backward switching on $G'$ is exactly the inverse of a forward switching which produces $G'$. Let $b(G')$ be the number of backward switchings that can be applied to $G'$. We will show the following.
\begin{claim}\label{c:backward}
\begin{eqnarray*}
(a)&&b(G')\le(d_u-d_u^{H_1})(d_v-d_v^{H_1})(M-2e(H_1)) \\
(b)&&b(G')\ge
(d_u-d_u^{H_1})(d_v-d_v^{H_1})(M-2e(H_1))\left(1- \frac{2J({\bf d})+6\Delta+2\Delta(H_2)\Delta}{M-2e(H_1)} \right).
\end{eqnarray*}
\end{claim}

Let $T$ be the total number of forward switchings from $S$ to $\overline{S}$. By definition,
$$
T = \sum_{G \in S} f(G) = \sum_{G' \in \overline{S}} b(G').
$$
\fix{Without loss of generality we may assume that $g({\bf d},H_1,H_2)>0$ as otherwise the lower bound on $\pr(uv\in \G(n,\bfd)\mid H_1^+, H_2^-)$ in Theorem~\ref{thm:conditional} holds trivially.}
By Claim~\ref{c:forward}(a) and Claim~\ref{c:backward}(b),
$$
|\overline{S}| (d_u-d_u^{H_1})(d_v-d_v^{H_1})(M-2e(H_1))\left(1- \frac{2J({\bf d})+6\Delta+2\Delta(H_2)\Delta}{M-2e(H_1)} \right) \le T \le |S|(M-2e(H_1))^2.
$$
Thus, 
\begin{eqnarray*}
\frac{|S|}{|S|+|\overline{S}|}&\ge& \frac{(d_v-d_v^{H_1})(d_u-d_u^{H_1})}{M-2e(H_1)} \left(1- \frac{2J({\bf d})+6\Delta+2\Delta(H_2)\Delta}{M-2e(H_1)} \right)
\left(1+ \frac{(d_v-d_v^{H_1})(d_u-d_u^{H_1})}{M-2e(H_1)}\right)^{-1}\\
&=&\frac{(d_v-d_v^{H_1})(d_u-d_u^{H_1})}{M-2e(H_1)}\cdot g(\bfd,H_1,H_2).
%\left[ \frac{1}{1-\left( \frac{2J({\bf d})+6\Delta+2\Delta(H^-)\Delta}{M-2e(H^+)} - \frac{(d_v-d_v^{H^+})(d_u-d_u^{H^+})}{M-2e(H^+)}\right)} \right].
\end{eqnarray*}
Similarly, by Claim~\ref{c:forward}(b) and Claim~\ref{c:backward}(a) \fix{and the assumption that $f({\bf d},H_1,H_2)>0$}, %and the fact that $(1-x)^{-1}\le 1+2x$ for all $x\in[0,1/2]$ yield
$$
\frac{|S|}{|S|+|\overline{S}|} \le \frac{(d_v-d_v^{H_1})(d_u-d_u^{H_1})}{M-2e(H_1)} \cdot f(\bfd,H_1,H_2).
$$
Hence, we have shown the upper and lower bounds of $\pr(uv\in \G(n,{\bf d}) \mid H_1^+, H_2^-)$ as desired.\qed \smallskip

It only remains to prove the two claims. They follow from simple inclusion-exclusion counting arguments as follows.\smallskip

{\em Proof of Claim~\ref{c:forward}. }
 The upper bound is obvious as there are at most $M-2e(H_1)$ ways to choose vertices $x$ and $a$, and then at most $M-2e(H_1)$ ways to choose vertices $y$ and $b$. To get the required lower bound, by the principle of inclusion and exclusion, it is sufficient to subtract from the above upper bound the number of choices where one of the conditions in (1)--(3) is violated. If condition (1) is violated, then $\{x, y, a, b\} \cap \{u, v\} \neq \emptyset$, or $a \in \{y, b\}$, or $x = b$. There are at most
    $2(M-2e(H_1))\cdot 2 \Delta \cdot 2 + 2\Delta(M-2e(H_1))+\Delta(M-2e(H_1))=11\Delta(M-2e(H_1))$ ways to choose such 4-tuples.
In our upper bound, we only considered choices where condition (2) is satisfied. Thus, it only remains to subtract the number of choices where condition (3) is violated. That means either (a),
 $xu$, $yv$, or $ab$ is an edge in $G$; or (b), $xu$, $yv$, or $ab$ is an edge in $H_2$. We call an ordered triple of vertices $(v_1,v_2,v_3)$ a {\em directed 2-path at $v_1$}, if both $v_1v_2$ and $v_2v_3$ are edges in the graph. Note that for any graph $G$ with degree sequence ${\bf d}$, and any $v\in[n]$, the number of directed 2-paths at $v$ in $G$ is always at most $\sum_{i=1}^{\Delta} (d_i-1)=J({\bf d})-\Delta$. Hence, the number of choices for (a) is at most
    $3(J({\bf d})-\Delta)(M-2e(H^+))$.
    The number of choices for (b) is at most
    $    (d^{H_2}_u\Delta+d^{H_2}_v\Delta)(M-2e(H_1))+2e(H_2)\Delta^2
$. These give the lower bound for $f(G)$ as desired.
\smallskip

{\em Proof of Claim~\ref{c:backward}. }
 The upper bound is obvious. There are at most $d_u-d_u^{H_1}$ ways to choose $x$, at most $d_v-d_v^{H_1}$ ways to choose $y$, and at most $M-2e(H_1)$ ways to choose $a$ and $b$. For the lower bound, by inclusion-exclusion,  we subtract from the above upper bound the number of choices where condition (1') or (2') is violated (note that our choices in the upper bound guarantee condition (3') already). If condition (1') is violated then $\{a, b\} \cap \{u, v, x, y\} \neq \emptyset$. There are at most
    $4\cdot 2\cdot(d_u-d_u^{H_1})(d_v-d_v^{H_1}) \Delta$ such choices.
If condition (2') is violated then either (a'), $xa$ or $yb$ is an edge in $G'$; or (b'), $xa$ or $yb$ is an edge in $H_2$. The number of choices for (a') is at most
    $2(d_u-d_u^{H_1})(d_v-d_v^{H_1})(J({\bf d})-\Delta)$,
    and the number of choice for (b')
    is at most
    $    2(d_u-d_u^{H_1})(d_v-d_v^{H_1})\Delta(H_2)\Delta$. Subtracting these upper bounds on the number of invalid choices from the upper bound on $b(G')$ yields the lower bound on $b(G')$ as desired.\qed

\section{Proof of Theorems~\ref{thm:small} and~\ref{thm:H}: connectivity}
\lab{sec:connectivity}

\subsection{Proof techniques: the old and the new}

Approximately four proof techniques or a hybrid of them have been used for proving the connectivity of $\G(n,\bfd)$ and for analysing properties of $\G(n,\bfd)$ in general, when $\Delta$ is not too large. The first, and perhaps the most well known method uses the configuration model. Recall that all a.a.s.\ results can be translated from the configuration model to $\G(n,\bfd)$ if $\bfd\in\D$ where $\D=\{\bfd: M=\Theta(n),\, \sum_{i=1}^n d_i^2=O(n)\}$. Federico and Hofstad's work~\cite{federico2016} is an example of this proof method. Due to the ease in handling with the configuration model, they managed to prove more accurate result including a critical window analysis during the connectivity phase transition, but such distributional results cannot be directly translated to $\G(n,\bfd)$.

Another proof method uses graph enumeration. Assume we want to bound the probability that $e(S,{\overline S})=0$. This probability is simply
\[
\frac{|\G(S,\bfd|_{S})| \cdot |\G({\overline S},\bfd|_{{\overline S}})|}{|\G(n,\bfd)|},
\] 
where $\bfd|_{S}$ denotes the degree sequence obtained by restricted to vertices in $S$, and $\G(S,\bfd|_{S})$ denote the set of graphs on vertex set $S$ and with degree sequence $\bfd|_{S}$. Applying known asymptotic enumeration results on $|\G(n,\bfd)|$ one can get asymptotic probability for the event that $e(S,{\overline S})$, which can further be used to bound the probability that $\G(n,\bfd)$ is disconnected. This approach was taken by Wormald~\cite{wormald1981}. It may be interesting to note that the enumeration results on which Wormald's proof was based are by Bender and Canfield~\cite{bender1978}, which requires $\Delta$ to be absolutely bounded. Then Bollob\'{a}s introduced the configuration model and deduced a probabilistic proof of~\cite{bender1978}. Afterwards the configuration model became popularised. %and many people analyse the configuration model directly instead of estimating edge probabilities using the enumeration results. 
In that sense, Wormald's proof can be viewed as a ``detour'' of the first method aforementioned.

The third method combines the configuration model with the switching method introduced by McKay~\cite{mckay1981,mckay1984}. As mentioned before, a.a.s.\ results can be translated from the configuration model to $\G(n,\bfd)$ if $\bfd\in \D$. What can we do for $\bfd\notin \D$, e.g.\ when the average degree of $\bfd$ is growing with $n$? McKay's switching method starts with $G^*$, the multigraph produced by the configuration model, and repeatedly switches away multiple edges in the multigraph from the configuration model. Using simple counting argument one can show that when $\Delta$ is below $M^{1/4}$ then the distribution of the final simple graph obtained is very close to the uniform distribution. Then, we can deduce properties of $\G(n,\bfd)$ by analysing the configuration model and the switching algorithm. There are many results of $\G(n,\bfd)$ obtained this way, e.g.\ $\chi(\G(n,\bfd))$ in~\cite{frieze2007} discussed in Section~\ref{sec:chrom}. See more examples in~\cite{janson2019}. In terms of the connectivity, proofs in \cite{frieze1988,luczak1989} followed this path.

The last method applies the switching technique directly to random graph $\G(n,\bfd)$. Partition the set of graphs $\G(n,\bfd)$ 
into two parts $S$ and $T$ where graphs in $S$ have a certain property $\P$ and graphs in $T$ do not. Then defining switchings that relate graphs in $S$ to graphs in $T$. By counting the number of ways to perform switchings one can estimate the ratio $|S|/|T|$ and  the probability of property $\P$. This approach was used in~\cite{cooper2002} for the connectivity of random $d$-regular graphs for $d$ up to $cn$ where $c$ is sufficiently small.

In this work, we use the new tool in Corollary~\ref{cor:upperbound} to characterise the connectivity phase transition for the family of degree sequences where $J(\bfd)=o(M)$ (Theorem~\ref{thm:small}). This result is a generalisation of~\cite{federico2016} but works for a much larger family of degree sequences. 
For Theorem~\ref{thm:H} we will use switchings to prove that the set of edges incident with $\cH$ spans a subgraph with $O(1)$ components. Then, we expose the set of edges incident with $\cH$, and then analyse the subgraph induced by $[n]\setminus\cH$. As the degrees of vertices in $[n]\setminus\cH$ are not too large, we can apply Corollary~\ref{cor:upperbound} again.

\subsection{Proof of Theorem~\ref{thm:small}(b,c)}

Let $Y$ denote the number of isolated edges, i.e.\ edges whose ends are both of degree 1. Let $Z$ be the number of isolated triangles, i.e.\ triples of vertices $\{x,y,z\}$ who induce a $K_3$ and all of the three vertices are of degree 2. With standard first and second moment calculations using the asymptotic probabilities in Corollary~\ref{cor:upperbound} we immediately have the following lemma, whose proof we omit.
\begin{lemma}
\begin{itemize}
\item If $n_1=\Omega(\sqrt{M})$ then $\ex Y\sim n_1^2/2M$ and $\ex Y(Y-1) \sim n_1^4/4M^2$.
\item If $n_2=\Omega(M)$ then $\ex Z\sim 4n_2^3/3M^3$ and $\ex Z(Z-1)\sim 16n_2^6/9M^6$.
\end{itemize}
\end{lemma}

Now Theorem~\ref{thm:small}(b) follows by Chebyshev's inequality, and Theorem~\ref{thm:small}(c) follows by the Paley-Zygmund inequality.\qed

\subsection{Proof of Theorem~\ref{thm:small}(a)}

\begin{prop}\lab{p:graphs}
For any $\bfd$ with even sum, 
\[
|\G(n,\bfd)|\le \frac{M!}{2^{M/2}(M/2)! \prod_{i\in [n]}d_i!}.  
\]
\end{prop}
\proof Represent vertex $i$ by a bin containing exactly $d_i$ points. A perfect matching over the total $M$ points in the $n$ bins is called a {\em pairing}. A pairing produces a multigraph with degree sequence $\bfd$ by representing each $\{v(p),v(q)\}$ as an edge where $p$ and $q$ are points matched by the pairing and $v(p)$ and $v(q)$ are the bins/vertices that contain points $p$ and $q$ respectively. It is easy to see that every graph in $\G(n,\bfd)$ corresponds to exactly $\prod_{i\in [n]} d_i!$ pairings. On the other hand, there are exactly $M!/2^{M/2}(M/2)!$ perfect matchings over $M$ points. The assertion follows.\qed\ss

Given $S\subseteq [n]$ let $X_S$ denote the indicator variable for the event that $e(S,{\overline S})=0$.   Recall that $d(S)=\sum_{i\in S} d_i$ and $\bfd|_S=(d_i)_{i\in S}$.

\begin{lemma}\lab{lem:AS} Assume $J(\bfd)=o(M)$.
Suppose $S\subseteq [n]$ where $M-d(S)=\Omega(M)$. Then, 
\[
\ex X_S\le \sqrt{2}\left(\frac{(1+o(1))d(S)}{M-d(S)}\right)^{d(S)/2}\left(1-\frac{d(S)}{M}\right)^{M/2}.
\]
\end{lemma}
\proof By Proposition~\ref{p:graphs}, the number of graphs on $S$ with degree sequence $\bfd|_S$ is at most
\[
\frac{d(S)!}{2^{d(S)/2}(d(S)/2)! \prod_{i\in S}d_i!}.  
\]
By~\eqn{graph_prob1},
\begin{eqnarray*}
\ex X_S &\le& \frac{d(S)!}{2^{d(S)/2}(d(S)/2)! \prod_{i\in S}d_i!} \cdot \prod_{i\in S} d_i! \prod_{i=1}^{d(S)/2}\frac{1+o(1)}{M-2i}\\
&\le&\frac{d(S)!}{2^{d(S)/2}(d(S)/2)!} \cdot (2+o(1))^{-d(S)/2} \frac{((M-d(S))/2)!}{(M/2)!} \\
&\le&\sqrt{2} \left(\frac{(1+o(1))d(S)}{M-d(S)}\right)^{d(S)/2}\left(1-\frac{d(S)}{M}\right)^{M/2},\qed
\end{eqnarray*}
\fix{where the last inequality holds by the Stirling approximation $\sqrt{2\pi n}\left(\frac{n}{e}\right)^n e^{\frac{1}{12n+1}} \le n! \le \sqrt{2\pi n}\left(\frac{n}{e}\right)^n e^{\frac{1}{12n}}$ for every $n\ge 1$.}

If $\G(n,\bfd)$ is disconnected then there is a component of $\G(n,\bfd)$ with total degree at most $M/2$. Thus, it is sufficient to show that a.a.s.\ there is no $S\subseteq [n]$ where $d(S)\le M/2$ and $e(S,{\overline S})=0$. In the next lemma, we first bound the expected number of such $S$ where $d(S)\ge 2.5|S|$.
\begin{lemma}\lab{lem:S}
A.a.s.\ there are no nonempty sets $S\subseteq [n]$ where $2.5|S|\le d(S)\le M/2$ and $e(S,{\overline S})=0$.
\end{lemma}

\proof Let $\eps=0.5$. Suppose $S$ is a set of vertices with $d(S)= h \le M/2$ and $d(S)\ge (2+\eps)|S|$. Then, $|S|\le h/(2+\eps)$. Hence, by Lemma~\ref{lem:AS}, for all sufficiently large $n$,

\begin{eqnarray*}
&&\sum_{\substack{S\subseteq [n]:\ 
d(S)= h\\
h\ge (2+\eps)|S|}}\ex X_S\le 2\left(\binom{n}{h/(2+\eps)}\ind{\frac{h}{2+\eps}<\frac{n}{2}} +2^n\ind{\frac{h}{2+\eps}\ge \frac{n}{2}} \right) \left(\frac{(1+o(1))\rho}{1-\rho}\right)^{h/2} (1-\rho)^{M/2},
\end{eqnarray*}
where $\rho=h/M< 1/2$.
By the assumption that $n_1=o(\sqrt{M})$ and $n_2=o(M)$, we must have $M\ge (3-o(1)) n$. Thus the above is at most
\begin{eqnarray*}
&&2\left((2e(1+\eps)/3)^{2/(2+\eps)} \frac{\rho^{\eps/(2+\eps)}}{1-\rho}\right)^{h/2}(1-\rho)^{M/2}+2(2^{2-\eps/2}\rho/(1-\rho))^{h/2}(1-\rho)^{M/2}.\\
%&&\hspace{0.5cm} \left(\left((2e(1+\eps)/3)^{2/(2+\eps)} \frac{\rho^{\eps/(2+\eps)}}{1-\rho}\right)^{\rho}(1-\rho)\right)^{M/2}+\Big((4\rho/(1-\rho))^{\rho}(1-\rho)\Big)^{M/2}
\end{eqnarray*}
We prove that
\[
\sum_{2\le h\le M/2}\left((2e(1+\eps)/3)^{2/(2+\eps)} \frac{\rho^{\eps/(2+\eps)}}{1-\rho}\right)^{h/2}(1-\rho)^{M/2}=o(1),
\]
and
\[
\sum_{2\le h\le M/2} (2^{2-\eps/2}\rho/(1-\rho))^{h/2}(1-\rho)^{M/2} =o(1),
\]
which will complete the proof of the lemma. Note that
\begin{eqnarray*}
&&\sum_{2\le h\le M/2}\left((2e(1+\eps)/3)^{2/(2+\eps)} \frac{\rho^{\eps/(2+\eps)}}{1-\rho}\right)^{h/2}(1-\rho)^{M/2}\\
&&=\sum_{2\le h<\ln n } \Big(O(1)(\ln n/M)^{\eps/(2+\eps)}\Big)^{h/2}+\sum_{\ln n\le h< 0.01 M}0.9^{h/2}+ \sum_{0.01 M\le h\le M/2}\exp\big(f(\rho)M/2\big),
\end{eqnarray*}
where
\[
f(\rho)=\frac{2\rho}{2+\eps}\ln(2e(1+\eps)/3) +\frac{\eps\rho}{2+\eps}\ln(\rho) +(1-\rho)\ln(1-\rho).
\]
The function $f(\rho)$ is below $-0.01$ uniformly over $\rho\in[0.01, 0.5]$. Hence,
\begin{eqnarray*}
\sum_{2\le h\le M/2}\left((2e(1+\eps)/3)^{2/(2+\eps)} \frac{\rho^{\eps/(2+\eps)}}{1-\rho}\right)^{h/2}(1-\rho)^{M/2}=o(1).
\end{eqnarray*}
Bounding $\sum_{2\le h\le M/2} (2^{2-\eps/2}\rho/(1-\rho))^{h/2}(1-\rho)^{M/2}$ by $o(1)$ can be done in a similar manner.\qed
\medskip

We are ready to complete the proof for part (a) of Theorem~\ref{thm:small}.

{\em Proof of Theorem~\ref{thm:small}(a). } By Lemma~\ref{lem:S} it only remains to show that a.a.s.\ there are no nonempty sets $S\subseteq [n]$ where $d(S)\le \min\{2.5|S|,M/2\}$ and $e(S,{\overline S})=0$. Let $\xi=\xi_n=o(1)$ be such that $n_1/\sqrt{M}\le \xi$ and $n_2/M\le \xi$. Let $S$ be a subset of vertices with $\ell_1$ vertices of degree 1, $\ell_2$ vertices of degree 2, and $\ell_{\ge 3}$ vertices of degree at least 3, $d(S)\le M/2$ and $d(S)\le 2.5|S|$. The number of ways to choose such a set $S$ is at most $\binom{n_1}{\ell_1}\binom{n_2}{\ell_2}\binom{n}{\ell_{\ge 3}}$. Given such an $S$, $d(S)\ge \ell_1+2\ell_2+3\ell_{\ge 3}$. It follows immediately that
\[
\ell_1+2\ell_2+3\ell_{\ge 3} \le 2.5(\ell_1+\ell_2+\ell_{\ge 3}),
\]
which implies that
\begin{equation}
\ell_{\ge 3}\le 3\ell_1+\ell_2, \quad \ell_1+2\ell_2+3\ell_{\ge 3}\le d(S)\le 2.5(\ell_1+\ell_2+\ell_{\ge 3})\le 10\ell_1+5\ell_2. \lab{h}
\end{equation}
In the rest of the proof, for simplicity we use $C$ to denote an absolute positive constant, which may take different values at different places where, the actual values of the constants do not matter. By Lemma~\ref{lem:AS}, the probability that $e(S,{\overline S})=0$ is at most
\[
\left(\frac{Cd(S)}{M}\right)^{d(S)/2} \le \left(\frac{C(\ell_1+\ell_2)}{M}\right)^{\ell_1/2+\ell_2+3\ell_{\ge 3}/2}.
\]
Hence,
  the expected number of sets $S$ where  $d(S)\le 2.5|S|$, $d(S)\le M/2$ and $e(S,{\overline S})=0$ is at most
\begin{eqnarray*}
&&\sum_{\ell_1,\ell_2,\ell_{\ge 3}} \binom{n_1}{\ell_1}\binom{n_2}{\ell_2}\binom{n}{\ell_{\ge 3}} \left(\frac{C(\ell_1+\ell_2)}{M}\right)^{\ell_1/2+\ell_2+3\ell_{\ge 3}/2}\\
&&\hspace{0.5cm} \le \sum_{\ell_1,\ell_2,\ell_{\ge 3}}\left(\frac{Cn_1}{\ell_1}\sqrt{\frac{\ell_1+\ell_2}{M}}\right)^{\ell_1} \left(\frac{Cn_2(\ell_1+\ell_2)}{\ell_2 M}\right)^{\ell_2} \left(\frac{Cn(\ell_1+\ell_2)^{3/2}}{\ell_{\ge 3} M^{3/2}}\right)^{\ell_{\ge 3}} \\
&&\hspace{0.5cm} \le \sum_{\ell_1,\ell_2,\ell_{\ge 3}}\left(C\xi \frac{\sqrt{\ell_1+\ell_2}}{\ell_1}\right)^{\ell_1} \left(C\xi \frac{\ell_1+\ell_2}{\ell_2}\right)^{\ell_2} \left(\frac{C(\ell_1+\ell_2)^{3/2}}{\ell_{\ge 3} \sqrt{M}}\right)^{\ell_{\ge 3}}. 
\end{eqnarray*}
We split the above sum into two parts, one restricted to $\ell_1\ge \ell_2$ and the other restricted to $\ell_1<\ell_2$, and we show that each sum is $o(1)$. Suppose $\ell_1 \ge \ell_2$. Then $\ell_1+\ell_2\le 2\ell_1$. Hence,
\begin{eqnarray*}
&&\sum_{\substack{\ell_1,\ell_2,\ell_{\ge 3}\\ \ell_1\ge \ell_2}}\left(C\xi \frac{\sqrt{\ell_1+\ell_2}}{\ell_1}\right)^{\ell_1} \left(C\xi \frac{\ell_1+\ell_2}{\ell_2}\right)^{\ell_2} \left(\frac{C(\ell_1+\ell_2)^{3/2}}{\ell_{\ge 3} \sqrt{M}}\right)^{\ell_{\ge 3}}\\
&&\hspace{0.5cm} \le \sum_{\substack{\ell_1,\ell_2,\ell_{\ge 3}\\ \ell_1\ge \ell_2}}\left(\frac{C\xi }{\sqrt{\ell_1}}\right)^{\ell_1} \left(C\xi \frac{\ell_1}{\ell_2}\right)^{\ell_2} \left(\frac{C\ell_1^{3/2}}{\ell_{\ge 3} \sqrt{M}}\right)^{\ell_{\ge 3}}.
\end{eqnarray*}
Let $g(x)=(K/x)^x$ on $x\ge 0$ where $K>0$. By considering the derivative of $\ln(g(x))$, it is easy to see that $g(x)$ is maximised at $x=K/e$.
Thus,
\[
\left(C\xi \frac{\ell_1}{\ell_2}\right)^{\ell_2}\le \exp\left(\frac{C\xi\ell_1}{e}\right), \quad\mbox{and}\quad  \left(\frac{C\ell_1^{3/2}}{\ell_{\ge 3} \sqrt{M}}\right)^{\ell_{\ge 3}} \le \exp\left(\frac{C\ell_1^{3/2}}{e\sqrt{M}}\right).
\]
It follows now that
\begin{eqnarray*}
&&\sum_{\substack{\ell_1,\ell_2,\ell_{\ge 3}\\ \ell_1\ge \ell_2}}\left(\frac{C\xi }{\sqrt{\ell_1}}\right)^{\ell_1} \left(C\xi \frac{\ell_1}{\ell_2}\right)^{\ell_2} \left(\frac{C\ell_1^{3/2}}{\ell_{\ge 3} \sqrt{M}}\right)^{\ell_{\ge 3}}= \sum_{\substack{\ell_1,\ell_2,\ell_{\ge 3}\\ \ell_1\ge \ell_2\\
\ell_1\le M^{1/3}/\log M}}\left(\frac{C\xi }{\sqrt{\ell_1}}\right)^{\ell_1} \left(C\xi \frac{\ell_1}{\ell_2}\right)^{\ell_2} \left(\frac{C\ell_1^{3/2}}{\ell_{\ge 3} \sqrt{M}}\right)^{\ell_{\ge 3}}\\
&&\hspace{1cm} 
+
\sum_{\substack{\ell_1,\ell_2,\ell_{\ge 3}\\ \ell_1\ge \ell_2\\
M^{1/3}/\log M<\ell_1\le n_1}}\left(\frac{C\xi }{\sqrt{\ell_1}}\right)^{\ell_1} \left(C\xi \frac{\ell_1}{\ell_2}\right)^{\ell_2} \left(\frac{C\ell_1^{3/2}}{\ell_{\ge 3} \sqrt{M}}\right)^{\ell_{\ge 3}}\\
&&\hspace{0.5cm}\le \sum_{
\ell_1\le M^{1/3}/\log M}\left(\frac{C\xi }{\sqrt{\ell_1}}\right)^{\ell_1} \cdot \ell_1 \exp\left(\frac{C\xi\ell_1}{e}\right)\cdot\sum_{\ell_{\ge 3}} \left(\frac{1}{\ell_{\ge 3} \log M}\right)^{\ell_{\ge 3}}\\
&&\hspace{1cm} 
+
\sum_{M^{1/3}/\log M<\ell_1\le n_1}\left(\frac{C\xi }{\sqrt{\ell_1}}\right)^{\ell_1}\cdot  \ell_1 \exp\left(\frac{C\xi\ell_1}{e}\right)  \cdot n\exp\left(\frac{C\ell_1^{3/2}}{e\sqrt{M}}\right)\\
&&\hspace{0.5cm}\le O\left(\frac{1}{\log M}\right) \sum_{
\ell_1\le M^{1/3}/\log M}\ell_1 \left(\frac{C\xi }{\sqrt{\ell_1}} e^{C\xi/e}\right)^{\ell_1} 
+
\sum_{M^{1/3}/\log M<\ell_1\le n_1}\ell_1 n \left(\frac{C\xi }{\sqrt{\ell_1}}e^{C\xi/e+C\sqrt{\ell_1/M}}\right)^{\ell_1} \\
&&\hspace{0.5cm} =o(1),
\end{eqnarray*}
as $n_1=o(\sqrt{M})$.

Similarly, we have
\begin{eqnarray*}
&&\sum_{\substack{\ell_1,\ell_2,\ell_{\ge 3}\\ \ell_1<\ell_2}}\left(C\xi \frac{\sqrt{\ell_1+\ell_2}}{\ell_1}\right)^{\ell_1} \left(C\xi \frac{\ell_1+\ell_2}{\ell_2}\right)^{\ell_2} \left(\frac{C(\ell_1+\ell_2)^{3/2}}{\ell_{\ge 3} \sqrt{M}}\right)^{\ell_{\ge 3}}\\
&&\hspace{0.5cm}\le \sum_{\substack{\ell_1,\ell_2,\ell_{\ge 3}\\ \ell_1<\ell_2}}\left(\frac{C\xi \sqrt{\ell_2}}{\ell_1}\right)^{\ell_1} \left(C\xi \right)^{\ell_2} \left(\frac{C\ell_2^{3/2}}{\ell_{\ge 3} \sqrt{M}}\right)^{\ell_{\ge 3}}\\
&&\hspace{0.5cm}\le \sum_{
\ell_2\le M^{1/3}/\log M} \ell_2 (C\xi e^{C\xi/\sqrt{\ell_2}} )^{\ell_2} \sum_{\ell_{\ge 3}}  \left(\frac{1}{\ell_{\ge 3} \log M}\right)^{\ell_{\ge 3}}\\
&&\hspace{1cm}+
 \sum_{
 M^{1/3}/\log M<\ell_2\le n_2} \ell_2n \big(C\xi e^{C\xi/\sqrt{\ell_2}+C\sqrt{\ell_2/M}} \big)^{\ell_2}\\
 &&\hspace{0.5cm}=o(1), 
\end{eqnarray*}
as $n_2=o(M)$. 
By Markov's inequality, a.a.s.\ there are no sets $S\subseteq [n]$ where $d(S)\le 2.5|S|$, $d(S)\le M/2$ and $e(S,{\overline S})=0$. This, together with Lemma~\ref{lem:S}, completes the proof for Theorem~\ref{thm:small}(a).
\qed

%%%%%%%%
%%%%%%%%

\section{Proof of Theorem~\ref{thm:H}}
Recall that
\[
\cH=\{i: d_i\ge \sqrt{M}/\log M\}.%,\quad {\widetilde M}=M-2d(\cH).
\]
We start by some structural result involving vertices in $\cH$.

\begin{lemma}\lab{lem1:connect} A.a.s.\ either all vertices in $\cH$ are contained in the same component of $\G(n,\bfd)$, or  $\G(n,\bfd)$ contains a set of components $\C_0,\C_1,\ldots,\C_k$ for some $2\le k\le \log^7 n$, where $\C_0$ contains no vertex in $\cH$, each $\C_j$ for $1\le j\le k$ contains at least one vertex in $\cH$, and $\cup_{j=1}^k \C_j$ has less than $M^{3/4}$ edges. 
\end{lemma}
\proof If $|\cH|\ge \log^7 M$ then all vertices in $\cH$ are contained in the same component of $\G(n,\bfd)$ by~\cite[Lemma 28]{joos2018}. Suppose $|\cH|<\log^7 M$. Let $\C_u$ denote the component of $\G(n,\bfd)$ that contains $u$. We prove that a.a.s.\ either all components in $\{\C_u, u\in\cH\}$ are of size smaller than $M^{2/3}$, or all vertices in $\cH$ are contained in the same component of $\G(n,\bfd)$. The proof is basically the same as in~\cite[Lemma 22]{joos2018}. We show that a.a.s.\ if there is $u\in \cH$ such that $|\C_u|\ge M^{2/3}$ then all vertices in $\cH$ must lie in the same component. Let $u,v\in \cH$. Let $\G^-$ be the set of graphs in $\G(n,\bfd)$ where $\C_u$ has at least  $M^{2/3}$ edges and $v\notin\C_u$. Let $\G^+$ be the set of graphs in $\G(n,\bfd)$ where $v\in\C_u$. We show that $\pr(\G^-)\le(\log M/M^{1/6})\pr(\G^+)$.

Define switchings between $\G^-$ and $\G^+$ as follows. Let $G\in \G^-$. Choose a neighbour $w$ of $v$ in $\C_v$. Choose two vertices $(x,y)$ in $\C_u$ such that $d(x,u)\le d(y,u)$ where $d(x,u)$ denotes the graph distance between $x$ and $u$ in $\C_u$. Delete the edges $vw$ and $xy$ and then add edges $xv$ and $wy$. Obviously the resulting graph $G'$ is in $\G(n,\bfd)$. Moreover, $u$ and $v$ lie in the same component. There are at least $\sqrt{M}/\log M$ ways to choose $w$ since $d_v\ge \sqrt{M}/\log M$. There are at least $M^{2/3}$ to choose $(x,y)$ since $\C_u$ has size at least $M^{2/3}$. Hence, the number of switchings is at least $M^{7/6}/\log M$. Let $G'\in \G^+$. Note that on all shortest paths in $G'$ between $u$ and $v$, the neighbours of $v$ must be the same, as otherwise $G'$ cannot be created by a switching from some $G\in \G^-$. Let $P$ be such a shortest path. Let $x$ be the neighbour of $v$ on $P$. Then there are at most $M$ ways to choose $(y,w)$ such that an inverse switching can be applied. Hence, the total number of inverse switchings is at most $M$. Hence,
\[
\frac{\pr(\G^-)}{\pr(\G^+)} \le \frac{M}{M^{7/6}/\log M} =\log M/M^{1/6}.
\]
Taking the union bound over all pairs of $u,v\in \cH$, it follows that a.a.s.\ if there is $u\in \cH$ such that $|\C_u|\ge M^{2/3}$ then all vertices in $\cH$ must lie in the same component. Suppose for all $u\in \cH$, $\C_u$ has less than $M^{2/3}$ edges then $\cup_{u\in\cH} \C_u$ has less than $M^{2/3}\log^7 n<M^{3/4}$ edges. It only remains to show that a.a.s.\ if $\cup_{u\in\cH} \C_u$ has less than $M^{3/4}$ edges then all vertices in $[n]\setminus (\cup_{u\in\cH} \C_u)$ lies in the same component $\C_0$.

Condition on the set of vertices $V'=[n]\setminus (\cup_{u\in\cH} \C_u)$ whose total degree $M'$ is at least $M-2M^{3/4}$. Moreover, all vertices in $V'$ have degree at most $\sqrt{M}/\log M=O(\sqrt{M'}/\log M')$. Let $n'_1$ and $n'_2$ be the number of vertices in $V'$ with degree one and two respectively. By the hypothesis of Theorem~\ref{thm:H}, $n_1'=o(\sqrt{M})=o(\sqrt{M'})$ and $n_2'=o(M)=o(M')$. By Theorem~\ref{thm:small}, a.a.s.\  all vertices in $V'$ are in the same component. \qed

\begin{lemma}\lab{lem2:connect} 
With probability $o(1)$, $\G(n,\bfd)$ contains a set of components $\C_0,\C_1,\ldots,\C_k$ for some $2\le k\le \log^7 n$, where $\C_0$ contains no vertex in $\cH$, each $\C_j$ for $1\le j\le k$ contains at least one vertex in $\cH$, and $\cup_{j=1}^k \C_j$ has less than $M^{3/4}$ edges. 
\end{lemma}
Lemmas~\ref{lem1:connect} and~\ref{lem2:connect} immediately give the following corollary.
\begin{cor}\lab{cor:connect}
A.a.s.\ all vertices in $\cH$ are contained in the same component of $\G(n,\bfd)$.
\end{cor}

{\em Proof of Lemma~\ref{lem2:connect}.\ } Let $\G^-$ be the set of graphs in $\G(n,\bfd)$ where there is a set of components $\C_0,\C_1,\ldots,\C_k$ for some $2\le k\le \log^7 n$, where $\C_0$ contains no vertex in $\cH$, each $\C_j$ for $1\le j\le k$ contains at least one vertex in $\cH$, and $\cup_{j=1}^k \C_j$ has less than $M^{3/4}$ edges.  Fix $u\in\cH$ and let $\G^+$ be the set of graphs in $\G(n,\bfd)$ where $u$ is in a component of size $\Omega(M)$. We show that
$\pr(\G^-)=o(\pr(\G^+))$ which then implies that $\pr(\G^-)=o(1)$. 

Let $G\in \G^-$. Let $\C_0,\ldots,\C_k$ be the set of components of $G$ where $\C_0$ contains no vertex in $\cH$, each $\C_j$ for $1\le j\le k$ contains at least one vertex in $\cH$, and $\cup_{j=1}^k \C_j$ has less than $M^{3/4}$ edges.  Define  switchings between $\G^-$ and $\G^+$ as follows. Let $v$ be a neighbour of $u$. Choose three vertices $(w,x,y)$ in $\C_0$ such that $xy$ in a non-bridge edge, and $d(x,u)\le d(y,u)$. Replace the edges $uv$ and $xy$ by $ux$ and $yv$. Obviously the resulting graph is in $\G^+$. Since $\C_0$ has size at least $M-M^{3/4}$, and the number of vertices in $\C_0$ of degree less than 3  is $o(M)$ by the hypothesis of Theorem~\ref{thm:H}, the number of choices for $(x,y)$ is $\Omega(M)$. The number of choices for $v$ is at least $\sqrt{M}/\log M$. The number of vertices in $\C_0$ is at least $n-2M^{3/4}$ since there are at most $M^{3/4}$ edges in components other than $\C_0$. The number of vertices in $\C_0$ is also bounded from below by $(M-2M^{3/4})\log M/\sqrt{M}\sim \sqrt{M}\log M$ since every vertex in $\C_0$ has degree at most $\sqrt{M}/\log M$. Hence, the number of choices for $w$ is at least $\varphi=\max\{n-2M^{3/4}, \frac{1}{2}\sqrt{M}\log M\}$. Thus, the number of switchings that can be applied to $G$ is $\Omega(M^{3/2}\varphi/\log M)$. For $G'\in\G^+$, we bound the number of inverse switchings that can be applied to $G'$. There are at most $n$ ways to choose $w$. After fixing $w$, let $P$ be a shortest path from $u$ to $w$ (as in the proof of Lemma~\ref{lem1:connect} we may assume on all such shortest paths the neighbours of $u$ are the same). Let $x$ be the neighbour of $u$ on $P$. Finally, there are at most $M$ ways to choose the edge $vy$ so that an inverse switching can be applied. Thus, the number of inverse switchings is at most $nM$. Hence,
\[
\frac{\pr(\G^-)}{\pr(\G^+)} =O\left(\frac{nM}{M^{3/2}\varphi/\log M}\right)= O\left(\frac{n\log M}{M^{1/2}\varphi}\right).
\]
Suppose $M\ge n^{4/3}/\log n$. Then $\varphi=\Omega(\sqrt{M}\log M)$ and thus the above ratio of probabilities is $O(n/M)=o(1)$.
Suppose $M<n^{4/3}/\log n$ then $\varphi\ge n-2M^{3/4}\sim n$ and thus the above ratio of probabilities is $O(\log M/M^{1/2})=o(1)$. It follows now that $\pr(\G^-)=o(1)$. \qed

\medskip

Now we are ready to prove the second main theorem of the paper.

{\em Proof of Theorem~\ref{thm:H}. } Expose the set of edges incident with $\cH$ and let ${\bf d}'=(d'_i)_{i\in[n]\setminus \cH}$ denote the remaining degree sequence for vertices in $[n]\setminus \cH$. That is, $d'_i=d_i-x_i$ where $x_i$ is the number of edges between $i$ and $\cH$. In this proof we will focus on the subgraph of $\G(n,\bfd)$ induced by $[n]\setminus \cH$. Conditioning on $\bfd'$, this subgraph is distributed as $\G([n]\setminus \cH,\bfd')$, a uniformly random graph on $[n]\setminus\cH$ with degree sequence $\bfd'$. Note that some vertices in $[n]\setminus \cH$ may have degree 0 with respect to $\bfd'$. They are not of interest for study as they are just isolated vertices in $\G([n]\setminus \cH,\bfd')$, and they are known to be adjacent to some vertex in $\cH$. Hence, let $V'$ be the set of vertices  $v\in [n]\setminus \cH$ where $d'_v\ge 1$.  Let $M'=\sum_{i\in V'}d'_i$. Conditioning on $V'$ and $\bfd'$, the subgraph of $\G(n,\bfd)$ induced by $V'$ is distributed as $\G(V',\bfd'|_{V'})$. By the theorem hypothesis that $M-2d(\cH)=\Omega(M)$, it follows that $M'\ge M-2d(\cH)=\Omega(M)$. Hence, by the definition of $\cH$ we have 
\begin{equation}
d'_i\le \sqrt{M}/\log M=O(\sqrt{M'}/\log M'),\quad\mbox{for all $i\in V'$.}  \lab{condition'}
\end{equation}
Given set $S\subseteq [n]$, we say $v$ is adjacent to $S$ if there exists $u\in S$ which is adjacent to $v$. 
\begin{claim}\lab{claim:H}
A.a.s.\  there exists $v\in [n]\setminus \cH$ such that $d_v\ge 2$ and $v$ is adjacent to $\cH$.
\end{claim}
 Let $U$ be the set of vertices $v\in [n]\setminus \cH$ where $v$ is adjacent to some vertex in $\cH$. Colour all vertices in $U$ red. The other vertices in $[n]\setminus \cH$ are uncoloured. Let $V_1=U\cap V'$. I.e.\ $V_1$ is the subset of vertices in $U$ with degree at least $1$ with respect to $\bfd'$. By Corollary~\ref{cor:connect}, we may assume that all red vertices are contained in the same component of $\G(n,\bfd)$. %By Claim~\ref{claim:H}, we may assume that $U\emptyset$. Then, the connectivity of $\G(n,\bfd)$ is implied if we can prove that there is no partition of $V'$ into $S$ and $T$ where $e(S,T)=0$, and no colour $i$ such that both $S\cap V_i\neq \emptyset$ and $T\cap V_i\neq\emptyset$. This implication is not true if there exists $i$ where $V_i= \emptyset$, as the set of vertices in $[n]\setminus \cH$ coloured $i$ (they are all isolated vertices with respect to $\bfd'$ since $V_i=\emptyset$) together with their neighbours in $\cH$ may lie in a distinct component from the vertices of other colours, and the uncoloured vertices. The next claim excludes such scenarios. 

\begin{claim}\lab{claim:connect0}
A.a.s.\  if there is a red vertex $u$ with $d'_u=0$ then $V_1\neq\emptyset$.
\end{claim}

Therefore we may assume that $V_1\neq \emptyset$.
In the rest of the proof, we will focus on $\G(V',\bfd'|_{V'})$, and we call $\bfd'_v$ the degree of $v$ for $v\in V'$. When we use graph notation such as $e(U,V)$ and $d(U)$, the graph referred to is $\G(V',\bfd'|_{V'})$ unless otherwise specified.  By construction, all vertices in $V'$ have degree at least 1. Moreover, by the theorem hypothesis on $n_1$ and $n_2$ and by the facts that $n_1'\le n_1$, $n_2'\le n_2$ and $M'=\Omega(M)$, where $n_1'$ is the number of uncoloured vertices of degree 1 in $V'$, and $n_2'$ is the number of uncoloured vertices of degree 2 in $V'$, it follows that 
\begin{equation}
n'_1=o(\sqrt{M'}),\quad n_2'=o(M'). \lab{n'} 
\end{equation}

Now, Theorem~\ref{thm:H} follows from the following  two claims.

\begin{claim}\lab{claim:connect1}
A.a.s.\ there is no $S\subset V'\setminus V_1$ where $d(S)\le M'/2$ and $e(S,V'\setminus S)=0$.
\end{claim}

\begin{claim}\lab{claim:connect2}
A.a.s.\ there exists no $T\subseteq [n]\setminus \cH$ where $V_1\subseteq T$, $d(T)\le M'/2$ and $e(T,V'\setminus T)=0$.

\end{claim}

The proofs of Claims~\ref{claim:connect1} and~\ref{claim:connect2} are analogous to the proof of Theorem~\ref{thm:small}. We briefly sketch the arguments.

{\em Proof of Claim~\ref{claim:connect1}. } By~\eqn{n'}, $n'_1\le \xi \sqrt{M'}$, and  $n'_2\le \xi M'$, for some $\xi=o(1)$. Moreover, by~\eqn{condition'}, the joint probabilities in~\eqn{graph_prob1} can be applied to $\G(V',\bfd'|_{V'})$. The rest of the proof is identical to that of Lemma~\ref{lem:S} and Theorem~\ref{thm:small}, noting that $S$ contains only uncoloured vertices. \qed
\medskip

{\em Proof of Claim~\ref{claim:connect2}. } Let $D$ denote the total degree of $V_{1}$, i.e.\ $D=\sum_{u\in V_1} d'_u$. Next, given a vector $\ell_1,\ell_2,\ldots$, the number of ways to choose $T$ where $V_1\subseteq T$, and there are $\ell_i$ uncoloured vertices of degree $i$ in $T$ is at most
\[
\binom{\xi \sqrt{M'}}{\ell_1} \binom{\xi M'}{\ell_2} \binom{n'}{\ell_{\ge 3}},
\]
where $\ell_{\ge 3}=\sum_{i\ge 3} \ell_i$ and $n'=|V'|$.
By Lemma~\ref{lem:S} we may assume that 
\[
d(T)=D+\sum_{i\ge 1} i\ell_i<2.5\left(\sum_{i\ge 1}\ell_i + D\right).
\]
The probability that $e(T,V'\setminus T)=0$ is at most
\begin{eqnarray*}
&&\left(\frac{d(T)}{M'-d(T)}\right)^{d(T)/2}\left(1-\frac{d(T)}{M'}\right)^{(M'-d(T))/2}\\
&&=\left(\frac{D+\sum_{i\ge 1} i\ell_i}{M'-(D+\sum_{i\ge 1} i\ell_i)}\right)^{(D+\sum_{i\ge 1} i\ell_i)/2}\left(1-\frac{D+\sum_{i\ge 1} i\ell_i}{M'}\right)^{(M'-(D+\sum_{i\ge 1} i\ell_i))/2}.
\end{eqnarray*}
Given $\ell_1,\ell_2,\ldots$, the above function is monotonely decreasing on $D$ on the domain where $(D+\sum_{i\ge 1}i\ell_i)/M'\le 1/2$. Hence, the above probability is maximised at $D=0$. Hence, the probability of existing such a set $T$, given $\ell_1, \ell_2,\ldots$, is at most
\begin{eqnarray*}
&&\binom{\xi \sqrt{M'}}{\ell_1} \binom{\xi M'}{\ell_2} \binom{n'}{\ell_{\ge 3}} \left(\frac{\sum_{i\ge 1} i\ell_i}{M'-\sum_{i\ge 1} i\ell_i}\right)^{\sum_{i\ge 1} i\ell_i/2}\left(1-\frac{\sum_{i\ge 1} i\ell_i}{M'}\right)^{(M'-\sum_{i\ge 1} i\ell_i)/2}\\
&&\le \binom{\xi \sqrt{M'}}{\ell_1} \binom{\xi M'}{\ell_2} \binom{n'}{\ell_{\ge 3}} \left(\frac{2\sum_{i\ge 1} i\ell_i}{M'}\right)^{\sum_{i\ge 1} i\ell_i/2}
\end{eqnarray*}
where $\sum_{i\ge 1} i\ell_i<2.5\sum_{i\ge 1}\ell_i $ implying 
\[
\ell_{\ge 3}\le 3\ell_1+\ell_2,\quad \ell_1+2\ell_2+3\ell_{\ge 3}\ell_3\le 2.5(\ell_1+\ell_2+\ell_{\ge 3})\le 10\ell_1+5\ell_2.
\]
The rest of the analysis is the same as in Theorem~\ref{thm:small}.\qed

\medskip

Now we prove Claims~\ref{claim:H} and~\ref{claim:connect0}. Both claims concern events related to edges incident with vertices in $\cH$. We will use switchings to bound probabilities of such events.

{\em Proof of Claim~\ref{claim:H}. } Let $\E$ denote the event that $e(\cH,[n]\setminus\cH)>0$ and all edges between $\cH$ and $[n]\setminus\cH$ have one end whose degree equals 1 in $\G(n,\bfd)$. If there is no vertex $v\in [n]\setminus\cH$ where $d_v\ge 2$ and $v$ is adjacent to $\cH$, then we must have either $\E$ or $e(\cH,[n]\setminus\cH)=0$. It is then sufficient to show that $\pr(\E)=o(1)$ and $\pr(e(\cH,[n]\setminus\cH)=0)=o(1)$. Let $\G$ be the class of graphs in $\G(n,\bfd)$ where $e(\cH,[n]\setminus\cH)=0$ and let $\G'$ be the class of graphs in $\G(n,\bfd)$ where $e(\cH,[n]\setminus\cH)=2$. If $\G=\emptyset$ then $\pr(e(\cH,[n]\setminus\cH)=0)=0$. Otherwise, 
\[
\pr(e(\cH,[n]\setminus\cH)=0) \le \frac{|\G|}{|\G'|}.
\]
We define a switching from $G\in \G$ by choosing an edge $xy$ in $G_{[\cH]}$ (note that $G_{[S]}$ denote the subgraph of $G$ induced by $S$) and another edge $uv$ in $G_{[[n]\setminus\cH]}$. Replace these two edges by $xu$ and $yv$. The resulting graph $G'$ is in $\G'$. There are $d(\cH)/2=\Omega(\sqrt{M}/\log M)$ ways to choose $xy$ and $d([n]\setminus \cH)=\Omega(M)$ ways to choose $uv$. So the number of switchings applicable on $G$ is at least $\Omega(M^{3/2}/\log M)$. On the other hand, for every $G'$, it can be produced by at most 1 way, as there are exactly 2 edges between $\cH$ and $[n]\setminus \cH$. It follows then that
\[
\pr(e(\cH,[n]\setminus\cH)=0)\le \frac{|\G|}{|\G'|}=O(\log M/M^{3/2}),
\]
as desired.

Next, let $\G$ be the class of  graphs in $\G(n,\bfd)\cap \E$, and let $\G'$ be the class of graphs in $\G(n,\bfd)$ where there is exactly one neighbour of $\cH$ in $[n]\setminus \cH$ with degree at least 2, and all the other neighbours of $\cH$ in $[n]\setminus \cH$ have degree equal to 1. Define a switching from $\G$ to $\G'$ as follows. Given $G\in \G$, choose 4 vertices $(u,v,x,y)$ such that $u\in\cH$,  $v,x,y\in [n]\setminus \cH$, $uv$ and $xy$ are edges, and $d_x\ge 2$. Since $d_x>1$ and $d_v=1$ it follows immediately that $ux$ and $vy$ are not edges. The switching replaces $uv$ and $xy$ by $ux$ and $vy$. The resulting graph $G'$ is in $\G'$ since $x$ becomes a neighbour of $\cH$ and $d_x\ge 2$. There is at least one way to choose $v$, since $G\in\E$. The total degree of $[n]\setminus \cH$ is $\Omega(M)$ by the theorem assumption, and there are at most $o(\sqrt{M})$ vertices in the set whose degree is 1. Moreover, all vertices in  $[n]\setminus \cH$ have degree at most $\sqrt{M}/\log M$. Hence, there are at least $\Omega(M)/(\sqrt{M}/\log M)=\Omega(\sqrt{M}\log M)$ vertices in $[n]\setminus \cH$ whose degree is at least 2. Hence, the number of choices for $x$ and $y$ is $\Omega(\sqrt{M}\log M)$. Thus, the number of switchings that can be applied to $G$ is $\Omega(\sqrt{M}\log M)$. On the other hand, given $G'\in\G$, there is a unique way to choose $u$ and $x$, and at most $n_1$ ways to choose $v$ and $y$ so that an inverse switching can be applied. Thus,
\[
\pr(\E)=O\left( \frac{n_1}{\sqrt{M}\log M}\right)=o(1),
\]
and our assertion follows.
\qed
\medskip

{\em Proof of Claim~\ref{claim:connect0}. }  Let $\P$ denote the set of graphs in $\G(n,\bfd)$ with the property in Claim~\ref{claim:H}, i.e.\ there exists $v\in[n]\setminus\cH$ where $d_v\ge 2$ and $v$ is adjacent to $\cH$. Let $\G$ denote the set of graphs in $\P$ where for all $v\in[n]\setminus \cH$ adjacent to $\cH$, $d'_v=0$. 
Let $\G'$ be the set of graphs in $\P$ where there are exactly two vertices $v_1,v_2\in[n]\setminus \cH$ such that $d'_{v_1}=1$, $d'_{v_2}\ge 1$, and $d'_v=0$ for all $v\in ([n]\setminus \cH)\setminus \{v_1,v_2\}$. Define a switching from $\G$ to $\G'$ as follows. Given $G\in\G$, the switching chooses 4 vertices $(u,v,x,y)$ such that $u\in \cH$, $v,x,y\in [n]\setminus \cH$, $uv$ and $xy$ are edges, $d_v\ge 2$ and $d'_x\ge 2$. Since $d'_x>0$ and $d'_v=0$, we know that $ux$ and $vy$ are not edges. The switching replaces edges $uv$ and $xy$ by $ux$ and $vy$. Let $G'$ denote the resulting graph. In $G'$, both $x$ and $v$ are adjacent to $\cH$, since $x$ is adjacent to $u\in\cH$, and $v$ is adjacent to some vertex $u'\in \cH$ where $u'\neq u$, as $d_v\ge 2$ and $d'_v$ was equal to 0 in $G$. Let $U'$ be the set of uncoloured vertices. The total degree of $U'$ is at least $M'=\Omega(M)$, since $e(U', [n]\setminus U')=0$ in $G$ by the definition of $\G$. Moreover, there are at most $o(\sqrt{M})$ vertices in $U'$ of degree one by~\eqn{n'}. Thus,  the number of choices for $x$ and $y$ is at least $\Omega(M)-o(\sqrt{M})=\Omega(M)$. Consequently, the total number of ways to perform a switching on $G$ is at least $\Omega(M)$. On the other hand, given $G'\in \G'$, there are exactly 2 vertices in $[n]\setminus \cH$ whose degrees (with respect to $\bfd'$) are at least 1, and they must be $v$ (the one with degree equal to 1) and $x$.  Fixing $v$ fixes $y$ as $d'_v=1$. Given $x$ there are at most $\sqrt{M}/\log M$ ways to choose $u$ as $d_x\le \sqrt{M}/\log M$. Hence, the number of ways $G'$ can be created via a switching is at most $2\sqrt{M}/\log M$. Together with Claim~\ref{claim:H}, 
\[
\pr(V_1=\emptyset) \le \frac{|\G|}{|\G'|}+\pr(\neg \P) \le \frac{2\sqrt{M}/\log M}{\Omega(M)}+o(1)=o(1).
\]
This completes the proof of the claim.\qed

%\bibliographystyle{plain}
%\bibliography{references}

\end{document}